\magnification 1200
\input amssym.def
\input amssym.tex
\parindent = 40 pt
\parskip = 12 pt
\font \heading = cmbx10 at 12 true pt
 at 22 true pt
\font \medheading =cmbx7 at 16 true pt
 at 7 true pt

\def \R{{\bf R}}

\centerline{\medheading Estimates for Fourier transforms of Surface  }
\centerline{\medheading Measures in R$^3$ and PDE applications }
\rm
\line{}
\line{}
\centerline{\heading Michael Greenblatt}
\line{}
\centerline{October 26, 2014}
\baselineskip = 12 pt
\font \heading = cmbx10 at 14 true pt
\line{}
\line{}
\line{}
\noindent{\heading 1. Background and Theorem Statements}.

\vfootnote{}{This research was supported in part by NSF grant  DMS-1001070}

Let $Q$ be a smooth two-dimensional surface in $\R^3$, and let $a$ be a point on $Q$. We consider the Fourier
transform of a small portion of the surface near $a$, localized using a smooth bump function supported near $a$. After a translation
and rotation, without loss of generality we may take $a = (0,0,0)$ and assume that $(0,0,1)$ is normal to $Q$ at the origin. In
this situation, we are looking at the following, where $\phi(x,y)$ denotes a smooth real-valued bump function supported near the origin and 
where $S(x,y)$ denotes the function whose graph is given by $Q$.
$$T(\lambda_1,\lambda_2,\lambda_3) = \int_{\R^2}e^{i\lambda_1S(x,y) +i\lambda_2x + i\lambda_3y}\,\phi(x,y)\,dx\,dy \eqno (1.1)$$
Technically this is the Fourier transform of the surface measure at $(-\lambda_1,-\lambda_2,-\lambda_3)$, but to simplify notation we will consider $T(\lambda_1,\lambda_2,\lambda_3)$ as written here. Note that $S(0,0) = 0$ and $\nabla S(0,0) = (0,0)$ due to our assumption that $(0,0,1)$ is normal to $Q$ at the origin. When $S(x,y)$ is flat to infinite order, one gets
very poor decay (if any) in $\lambda_1$ when $\lambda_2 = \lambda_3 = 0$ and there can be other pathologies, so we always
assume that at least one partial $\partial_x^{\alpha}\partial_y^{\beta}S(0,0) \neq 0$.

When $\lambda_2 = \lambda_3 = 0$, the function $U(\lambda_1) = T(\lambda_1,0,0)$ becomes a standard scalar oscillatory integral,
and it is well-known (see [AGV] ch. 6) that when $S(x,y)$ is real-analytic, if $\phi$ is supported in a sufficiently small neighborhood of the origin then as $\lambda_1 \rightarrow \infty$ one has an asymptotic development of the form
$$U(\lambda_1) = c_{S,\phi} \lambda_1^{-{\epsilon}}(\ln \lambda_1)^m + o(\lambda_1^{-{\epsilon}}(\ln \lambda_1)^m) \eqno (1.2)$$
Here $m = 0$ or $1$, and the pair $(\epsilon, m)$ is independent of $\phi$ and determined by the resolution of singularities
of $S(x,y)$ at the origin. The constant $c_{S,\phi}$ will be nonzero whenever $\phi$ is nonnegative with $\phi(0,0) > 0$. When
$\lambda_1$ is negative, then $U(\lambda_1)$ is just the complex conjugate of $U(|\lambda_1|)$ . Thus an expansion of the
form $(1.2)$ in $|\lambda_1|$ also holds as $\lambda_1 \rightarrow -\infty$.

In the more general smooth case, in [G1] it is shown
there is always an $(\epsilon, m)$ with $\epsilon > 0$ and $m = 0$ or $1$ such that for $|\lambda_1| > 2$  one has an upper bound
$$|U(\lambda_1)| \leq  c_{S,\phi} |\lambda_1|^{-{\epsilon}}(\ln|\lambda_1|)^m  \eqno (1.3)$$
We stipulate that $|\lambda_1| > 2$ to avoid trivial cases where one has to change the formula due to the fact that $\ln(1) = 0$.
This $(\epsilon, m)$ has the property that if $\phi$ is nonnegative with $\phi(0,0) > 0$, then $(1.3)$ will not hold for $(\epsilon',m')$ with $\epsilon' > \epsilon$, or for $\epsilon' = \epsilon$ with $m' < m$. It was then shown in [IM2] that most of the time one even has a development of the form $(1.2)$. 
It was shown in [V] for the real-analytic case and in [IM1] for the general smooth case that there are always certain "adapted" coordinate systems in which one can read off $(\epsilon,m)$ in terms of the Newton polygon of $S(x,y)$, and criteria can be given to determine if one is in such an adapted coordinate systems.

Note that in $(1.1)$, if  for a given $\delta > 0$ one has $|\lambda_2| + |\lambda_3| > \delta |\lambda_1|$ then if  the support of $\phi$ is sufficiently small (depending on $\delta$) the gradient
of the phase in $(1.1)$ is nonvanishing throughout the support of $\phi$. Thus one can do repeated integrations by parts and
for any $N$ one can quickly get an estimate of the form $|T(\lambda)| < C_N|\lambda|^{-N}$. Thus one always assumes that
$|\lambda_2| + |\lambda_3| \leq \delta |\lambda_1|$ for some small but fixed $\delta$. In part because of this, in much of the work concerning the oscillatory integrals $(1.1)$, people have viewed $(1.1)$ as perturbations of $U(\lambda_1)$ and proven upper bounds of the form $|T(\lambda)| \leq
C_{S, \phi} |\lambda_1|^{-{\epsilon}}(\ln|\lambda_1|)^m$, where $\epsilon$ and $m$ are as in $(1.2)$ or $(1.3)$. In particular, in the real-analytic case, a theorem of 
Karpushkin [K1]-[K2] says that one always has upper bounds of this form. In the smooth situation, for the case where $\epsilon > {1 \over 2}$ such upper bounds are a consequence of [D], and for the $\epsilon \leq {1 \over 2}$ situation these upper bounds are proven in [IKeM] [IM2]. One can obtain stronger results if one restricts to specific classes of functions, such as when the Hessian determinant is nonzero (where one has the strongest decay), the convex case considered in [BNW] [CoMa], or the class of surfaces in [ESa]. Curvature has often played a prominent role in such theorems. Other oscillatory integrals related to surface measure Fourier transforms were analyzed in [Gr].

We now let $\mu = (\lambda_2,\lambda_3)$, so that $\lambda$ may be written as $(\lambda_1,\mu)$. Our first theorem says that
in the general real-analytic case, one has $|T(\lambda_1,\mu_1,\mu_2)|  < C_{S,\phi}|\mu|^{-{1 \over 2}}$. It goes beyond what follows
from the perturbation results (Karpushkin's theorem) when $\epsilon \leq {1 \over 2}$.

\noindent {\bf Theorem 1.1.} Suppose $S(x,y)$ is real-analytic. There is a neighborhood $V$ of the origin such that if $\phi$ is
supported in $V$ then for some constant $C_S$ we have the following, where $|\mu|$ denotes
the magnitude of the vector $(\mu_1,\mu_2)$.
$$|T(\lambda_1,\mu_1,\mu_2)| < C_S |\mu|^{-{1 \over 2}}\|\phi\|_{C^1(V)} \eqno (1.4)$$

It can be shown that for many specific phases one gets a better exponent than ${1 \over 2}$ in $(1.4)$, but ${1 \over 2}$ is
the best exponent that holds for all phases, as can be seen when $S(x,y)$ is a function of $x$ or $y$ only.
Typically one does not expect to get a better exponent than $1$. This is because that is the decay rate for nondegenerate phases, so if one chooses $\phi$ supported in a small ball where $\nabla S$ and the Hessian determinant of $S$ are nonvanishing (assuming one exists) one will get a decay rate $\sim |\mu|^{-1}$, which can be seen by examining the $|\lambda_1| \sim |\mu|$ range and letting $|\mu|, |\lambda_1| \rightarrow \infty$.  

The next theorem will provide a new proof of the perturbation results for general smooth phase when $\epsilon \leq {1 \over 3}$. In the terminology of
 Varchenko [V] and later papers, this corresponds to
when the height of $S$ is at least 3.  Although such results are known in the real-analytic case by [K1][K2], and in the general
smooth case by [IKeM][IM2], we give a new proof here to 
illustrate that such theorems can also be proven with an appropriate
resolution of singularities theorem, without reference to adapted coordinates and so on. While
there are certainly commonalities between the proof of Theorem 1.2 and the arguments in [IKeM][IM2], there are also noteworthy
differences due to the use here of the resolution of singularities theorem of the next section and its consequences such as Lemma 2.2, as opposed to the type of subdivisions made in those papers.

\noindent {\bf Theorem 1.2.} Suppose $S(x,y)$ is smooth and $(\epsilon,m)$ is as in $(1.3)$. If $\epsilon \leq {1 \over 3}$, then there is a neighborhood $V$ of the origin such that if $\phi$ is supported in $V$ then for $|\lambda_1| > 2$ one has 
$$|T(\lambda_1,\mu_1,\mu_2)| \leq C_S |\lambda_1|^{-{\epsilon}}(\ln|\lambda_1|)^m  \|\phi\|_{C^1(V)} \eqno (1.5)$$
Again the $|\lambda_1| > 2$ condition is here to avoid concerning ourselves with trivial cases where one has to change the formula due to the fact that $\ln(1) = 0$.

\noindent {\heading PDE applications.}

We now assume $S(x_1,x_2)$ is real-analytic on some open ball $B$ centered at the origin with $S(0,0) = 0$ and $\nabla S(0,0) = (0,0)$. Suppose $f(x)$ is a complex-valued function on $\R^2$ such that $\hat{f}(\xi)$ is $L^1$ and is supported in $B$.  Let $F$ denote the Fourier transform, and define $S(-i \partial)$ to be the operator such that $F(S(-i \partial) f)(\xi) = S(\xi)\hat{f}(\xi)$. When $f$ is a function of $(t,x_1,x_2)$ we interpret this to be this multiplier operator in the 
$x_1$ and $x_2$ variables, with $t$ fixed. In Section 5, using Theorem 1.1 along with the general real-analytic version of Theorem 1.2  (i.e. Karpushkin's work) we will prove the following.

\noindent {\bf Theorem 1.3.} Suppose $(\epsilon,m)$ is as in $(1.3)$ for a real-analytic $S(x,y)$ and $\epsilon \leq {1 \over 2}$. If $B$ is sufficiently small, then the following holds. Let $1 \leq p  < \infty$. For $g$ such that $\hat{g} \in C_c^{\infty}(B)$, let $f(t,x_1,x_2)$ be the 
solution on $\R^3$ to the partial differential equation
$${\partial f \over \partial t}(t,x_1,x_2) = iS(-i\partial) f(t,x_1,x_2)$$
$$f(0,x_1,x_2) = g(x_1,x_2) \eqno (1.6)$$
Then  if $1 < q \leq \infty$ satisfies ${1 \over q} - {1 \over p} + {3 \over 4} < 0$ there is a constant $C_{p,q,S}$ such that one has the estimate
$$ \|f\|_q \leq C_{p,q,S}(|t| + 2)^{4\epsilon({1 \over q} - {1 \over p} + {3 \over 4})}(\ln(|t| + 2))^{-4m({1 \over q} - {1 \over p} + {3 \over 4})}\|g\|_p
\eqno (1.7)$$
The same is true if ${1 \over q} - {1 \over p} + {3 \over 4} = 0$, as long as $p \neq 1$ and $q \neq \infty$. Here the $L^p$ and $L^q$ norms are in the $x$ variables.

We have the condition $\epsilon \leq {1 \over 2}$ in Theorem 1.3 since when $\epsilon > {1 \over 2}$ one can get a stronger result by relatively rudimentary means. Consider now the case where $S(x_1,x_2) \geq 0$ in a neighborhood of the origin, and consider the oscillatory integral $R(\lambda_1, \mu_1, \mu_2)$ defined by
$$ R(\lambda_1, \mu_1, \mu_2)  = \int_{\R^2}e^{-\lambda_1S(x_1,x_2) +i\mu_1 x_1 + i\mu_2x_2}\,\phi(x_1,x_2)\,dx_1\,dx_2 \eqno (1.8)$$
In other words, we replace the $i\lambda_1S(x_1,x_2)$ in $(1.1)$ by $-\lambda_1S(x_1,x_2)$. In Lemma 5.1, we will see (by a much
easier argument than those proving Theorems 1.1 and 1.2) that if the support of $\phi$ is sufficiently small then for $\lambda_1 > 2$ one has an estimate
$$ |R(\lambda_1, \mu_1, \mu_2)| \leq C_{S,\phi} \min(\lambda_1^{-\epsilon} (\ln \lambda_1)^m, |\mu|^{-1}) \eqno (1.9)$$ 
Then in analogy to Theorem 1.3, one has the following theorem.

\noindent {\bf Theorem 1.4.} Assume $S$ is real-analytic and nonnegative on a neighborhood of the origin.  If $B$ is sufficiently small, then the following holds. Let $1 \leq p < \infty$. For $g$ such that $\hat{g} \in C_c^{\infty}(B)$,  let $f(t,x_1,x_2)$ be the 
solution on $\R^3$ to the partial differential equation
$${\partial f \over \partial t}(t,x_1,x_2) = -S(-i\partial) f(t,x_1,x_2)$$
$$f(0,x_1,x_2) = g(x_1,x_2) \eqno (1.10)$$
Then if $1 < q \leq \infty$ satisfies ${1 \over q} - {1 \over p} + {1 \over 2} < 0$, there exists a constant $C_{p,q,S}$ such that for $t > 0$  one has estimate
$$ \|f\|_q \leq C_{p,q,S}(t + 2)^{2\epsilon({1 \over q} - {1 \over p} + {1 \over 2})}(\ln(t + 2))^{-2m({1 \over q} - {1 \over p} + {1 \over 2})}\|g\|_p
\eqno (1.11)$$
The same is true if ${1 \over q} - {1 \over p} + {1 \over 2} = 0$, as long as $p \neq 1$ and $q \neq \infty$.
Here $(\epsilon,m)$ is as in $(4.2)$ for $S(x,y)$, and the $L^p$ and $L^q$ norms are in the $x$ variables.

Note that in Theorem 1.4, $(\epsilon, m)$ is as in $(4.2)$ and not as in $(1.3)$. By [G1] these may differ only when $(\epsilon,m) = (1,0)$ in $(1.3)$;  in this case the $(\epsilon, m)$ in $(4.2)$ may be either $(1,0)$ or $(1,1)$.

Theorem 1.4 can be used to relatively quickly give the following. Here $[S(-i\partial)]^{\delta}$ refers to the operator with
Fourier multiplier $(S(\xi_1,\xi_2))^{\delta}$ (we will only be considering it on $\xi$ domains where $S(\xi_1,\xi_2)$ is nonnegative).

\noindent {\bf Theorem 1.5.} Again assume $S$ is real-analytic and nonnegative on a neighborhood of the origin, and again let $(\epsilon,m)$ be as in $(4.2)$ for $S(x,y)$. If $B$ is sufficiently small, then the following holds. Let $0 < \delta < \epsilon$. For $g$ such that $\hat{g} \in C_c^{\infty}(B)$, let $f(x_1,x_2)$ solve the equation 
$$[S(-i\partial)]^{\delta}f = g \eqno (1.12)$$ 
Then if $p\in [1,\infty)$ and $q \in (1,\infty]$ such that  ${1 \over q} - {1 \over p} + {1\over 2} +  {\delta \over 2 \epsilon} < 0$, one has an estimate of the form
$$\|f\|_q \leq C_{p,q,S}\|g\|_p \eqno (1.13)$$
When $m = 0$, the same is true if  ${1 \over q} - {1 \over p} + {1\over 2} +  {\delta \over 2 \epsilon} = 0$, so long as $p \neq 1$ and
$q \neq \infty$.

The condition that  $\delta < \epsilon$ is needed in Theorem 1.5 for the statement to make sense; if $\delta \geq \epsilon$ then
$S(\xi_1,\xi_2)^{-\delta}$ is not integrable on a neighborhood of the origin and one cannot even automatically refer to the solution to $(1.12)$.

Theorems 1.3-1.5 are not intended to give the best possible exponents, or in the case of Theorems 1.3-1.4, the best 
possible powers of $|t|$ and $\ln |t|$, for any particular $S(x,y)$. Rather, they are illustrations of how one may interpret in terms of PDE theorems the combination of Theorem 1.1 and 1.2 or their analogues for $R(\lambda_1, \mu_1, \mu_2)$, in such a way as to
give estimates for any given $S(x,y)$.

\noindent {\heading 2. The Resolution of Singularities Theorem.} 

We next describe the resolution of singularities theorem that we need for this paper. There have been various resolution of 
singularities algorithms used in classical analysis problems in two dimensions, such as those of [PS] [V] and the author's 
earlier work. For the purposes of this paper we will use a modification of the one used in [G2], which was influenced by 
both [PS] and [V].

Suppose $f(x,y)$ is any smooth function on a neighborhood of the origin such that $f(0,0) = 0$ and such
that the Taylor expansion of $f$ at the origin has at least one nonvanishing term. After a linear change of coordinates if necessary
we may assume that if $k$ denotes the order of the zero of $f(x,y)$ at the origin then the Taylor expansion $\sum_{\alpha \beta}f_{\alpha\beta}x^{\alpha}y^{\beta}$ of $f$ at the origin
contains both a nonvanishing $f_{k0}x^k$ term and a nonvanishing $f_{0k}y^k$ term. We will now apply the resolution of singularities algorithm of Theorem 3.1 of [G2] in the following fashion. We divide the
$xy$ plane into 8 triangles via the $x$ and $y$ axes as well as two lines through the origin, one of the form $y = mx$ for $m > 0$ and one of the form $y = mx$ for $m < 0$. For certain technical reasons, these two lines cannot be ones on 
which the function $f_0(x,y) = \sum_{\alpha + \beta = k}f_{\alpha\beta}x^{\alpha}y^{\beta}$ vanishes. After 
possible reflections about the $x$ and/or $y$ axes and/or the line $y = x$, modulo its boundary each of the triangles is of the form $B_a = \{(x,y) \in \R^2: x > 0,\,\, 0 < y < ax\}$.

We now apply Theorem 3.1 of [G2] to the (reflected) $f(x,y)$ on the portion of $B_a$ contained in a sufficiently small neighborhood of
the origin. Actually, we apply a slight variant. If in the first step of the proof of Theorem 3.1 of [G2] one does a coordinate change of 
the form $(x,y) \rightarrow (x, y + cx + $ higher order terms$)$, instead we just do a coordinate change $(x,y) \rightarrow 
(x,y + cx)$. This has some technical advantages; see the proof of Theorem 2.1 d) below. Other than this, we do exactly the algorithm of Theorem 3.1  of [G2]. The following theorem is then a consequence of Theorem 3.1 of [G2].

\noindent {\bf Theorem 2.1.} Suppose $B_a  = \{(x,y) \in \R^2: x > 0,\, 0 < y < ax\}$ is as above. Abusing notation slightly, use the notation $f(x,y)$ to denote the reflected function $f(\pm x,\pm y)$ or $f(\pm y, \pm x)$ corresponding to $B_a$.
 Then there is a $b > 0$ and a positive integer $N$ such that
if $F_a$ denotes  $\{(x,y) \in \R^2: 0 \leq x\leq b, \,0 \leq y \leq ax\}$, then one can write $F_a = \cup_{i=1}^n cl(D_i)$, such that for to each $i$ there is a $\psi_i(x)$ with $\psi_i(x^N)$ smooth and $\psi_i(0) = 0$ such that after a coordinate change of the form $\eta_i(x,y) = (x, \pm y + \psi_i(x))$, the set $D_i$ becomes a set $D_i'$ on which the function $f \circ \eta_i(x,y)$ approximately becomes a monomial $d_i x^{\alpha_i}y^{\beta_i}$, $\alpha_i$ a nonnegative rational number and $\beta_i$ a nonnegative integer as follows.

\noindent {\bf a)} $D_i' = \{(x,y): 0 < x < b, \, g_i(x) < y < G_i(x)\}$, where $g_i(x^N)$ and $G_i(x^N)$ are
smooth. If we expand $G_i(x) =  H_i x^{M_i} + ...$, then $M_i \geq 1$ and $H_i > 0$, and consists of a single term $H_ix^{M_i}$ when $\beta_i = 0$. The function $g_i(x)$ is either identically
zero or can be expanded as $h_ix^{m_i} + ...$ where $h_i > 0$ and $m_i > M_i$.

\noindent {\bf b)} If $\beta_i = 0$, then $g_i(x)$ is identically zero. Furthermore, the $D_i'$ can
be constructed such that for any predetermined $\delta > 0$ there is a $d_i \neq 0$ such that on $D_i'$, for all $0 \leq l \leq \alpha_i$ one has
$$ |\partial_x^l (f \circ \eta_i)(x,y) -  d_i\alpha_i(\alpha_i -1) ... (\alpha_i - l + 1)x^{\alpha_i - l}| < \delta|d_i| x^{\alpha_i-l} \eqno (2.1)$$
This $\delta$ can be chosen independent of all the exponents appearing in this theorem. Furthermore, if one Taylor expands $f \circ \eta_i(x,y)$ in powers of $x^{1 \over N}$ and $y$ as  $\sum_{\alpha,\beta}F_{\alpha,\beta}x^{\alpha}y^{\beta}$, then $\alpha_i \leq \alpha + M_i\beta$ for all $(\alpha,\beta)$ such that $F_{\alpha,\beta} \neq 0$, with equality holding for at least two
$(\alpha,\beta)$, one of which is $(\alpha_i,0)$ and another of which satisfies $\beta > 0$. 

\noindent {\bf c)} If $\beta_i > 0$, then one may write $f = f_1^i + f_2^i$ as follows. $f_2^i \circ \eta_i(x,y)$ has a zero of infinite order at $(0,0)$ and
is identically zero if $f$ is real-analytic. $f_1^i \circ \eta_i(x^N,y)$ is smooth and there exists a $d_i \neq 0$ such that for any predetermined $\delta > 0$ (which can be chosen independent of the exponents appearing in this theorem) the $D_i'$ can
be constructed such that on $D_i'$, for any $0 \leq l \leq \alpha_i$ and any $0 \leq  m \leq \beta_i$ one has
$$|\partial_x^l\partial_y^m(f_1^i \circ \eta_i)(x,y) -  \alpha_i(\alpha_i - 1)....(\alpha_i - l + 1)\beta_i(\beta_i - 1)...(\beta_i - m + 1)
d_ix^{\alpha_i - l}y^{\beta_i - m}| $$
$$\leq \delta |d_i| x^{\alpha_i - l}y^{\beta_i - m} \eqno (2.2)$$

\noindent {\bf d)} If $\beta_i = 0$ and we write $\psi_i(x) = k_ix^{r_i} + ...$, then either $\psi_i(x) = k_ix$ for some $k_i$, $\psi_i(x) = k_ix + l_ix^{s_i}$ with $s_i = M_i > 1$ and $l_i \neq 0$,  or $\psi_i(x) = k_i x + l_ix^{s_i} + $ higher-order terms (if any), where $l_i\neq 0$ and $M_i > s_i > 1$. 

\noindent {\bf Proof.} Part a) is part of the statement of Theorem 3.1 of [G2], other than the form of the upper
edge of $D_i'$ when $\beta_i = 0$, which is given in the proof itself. Part c) is also contained in the statement of Theorem 3.1 of [G2]. 

As for part b), a weaker version was proved in [G2] using equation $(3.4)$ of that paper, and the stronger statement here also follows from that equation; if one divides $D_i'$ into finitely many subwedges of
width $\sim \epsilon x^{M_i}$ for small $\epsilon$ and then does a coordinate change of the form $(x,y - cx^{M_i})$ on each subwedge that transfers
its lower boundary to the $x$-axis, then if the subwedges are narrow enough, equation $(3.4)$ of [G2] implies that $(2.1)$ holds. Decreasing $\epsilon$ ensures that $\delta$ can be made as small as one would like.
As for the last sentence of
part b), although it is not in the statement of Theorem 3.1 of [G2] it is shown in the proof. 

Part d) is a consequence of the fact that in the version of the algorithm here, for a $D_i'$ with $\beta_i = 0$ one starts with a coordinate change of the form $(x,y) \rightarrow (x,y + k_ix)$, $k_i \neq 0$ if needed.  If additional coordinate changes are needed, then the 
 second coordinate change is either of the form $(x, \pm y + l_ix^{M_i})$ with $M_i = s_i$ and we are done, or it is of the form $(x, \pm y + l_ix^{s_i} + $ possible higher order terms) in such a way that the domains eventually giving a $\beta_i = 0$ wedge already are of width
$cx^m$ for some $m > s_i$. Further iterations of the resolution of singularities process will only add terms
of degree greater than $s_i$ and narrow the wedge further, resulting in an $M_i > s_i$.

\noindent The next lemma is a consequence of Theorem 2.1 we will need for our arguments.

\noindent {\bf Lemma 2.2.} Suppose $D_i'$ is such that $\beta_i = 0$. Then on $[0,b] \times [0,H_i]$ we may write
$$f \circ \eta_i(x, x^{M_i}y) = x^{\alpha_i}r_i(y) + E_i(x,y) \eqno (2.3)$$
Here $r(y)$ is a polynomial that doesn't vanish on $[0,H_i]$ and there is a $\delta > 0$ such that for any $l \geq 0$ there is a constant $C_{il}$ such that $E_i(x,y)$ satisfies
$$|\partial_x^l E_i(x,y)| \leq C_{il}x^{\alpha_i + \delta - l} \eqno (2.4)$$
\noindent {\bf Proof.} Again write $f \circ \eta_i (x,y) = \sum_{\alpha,\beta} F_{\alpha,\beta}x^{\alpha}y^{\beta}$. By part 
b) of Theorem 2.1, the minimum of $\alpha + M_i\beta$ in the sum above is $\alpha_i$ and furthermore $F_{\alpha_i,0} \neq 0$.
Let $q_i(x,y)$ be the polynomial $\sum_{\alpha + M_i\beta = \alpha_i} F_{\alpha,\beta}x^{\alpha}y^{\beta}$. Then by mixed
homogeneity we may write $q_i(x,y) = x^{\alpha_i}q_i(1, {y \over x^{M_i}})$. We now do a partial Taylor expansion of $f \circ \eta_i (x,y)$ in the form
$$f \circ \eta_i (x,y) = q_i(x,y) + \sum_{\alpha_i < \alpha + M_i\beta < K} F_{\alpha,\beta}x^{\alpha}y^{\beta} + O(x^K) 
\eqno (2.5)$$
Here $K$ is a large number determined by our arguments. We have an $O(x^K)$  and not an $O(x^K) + O(y^K)$ remainder term here because $0 < y <  H_ix^{M_i}$ on $D_i'$. Next, note that $(2.5)$ implies
$$f \circ \eta_i (x,x^{M_i}y) = x^{\alpha_i}q_i(1,y)+ \sum_{\alpha_i < \alpha + M_i\beta < K}F_{\alpha,\beta}x^{\alpha + M_i\beta}y^{\beta} + O(x^K) \eqno (2.6)$$
By Theorem 2.1 b), there are positive constants $e_i$ and $E_i$  such that on $(0,b] \times [0,H_i]$ one has
$$ e_i \leq {|f \circ \eta_i (x,x^{M_i}y)| \over x^{\alpha_i}}  \leq  E_i \eqno (2.7)$$
So by dividing by $x^{\alpha_i}$ and taking limits as $x \rightarrow 0$ in $(2.6)$ we have that $q_i(1,y) \neq 0$ for $0 \leq  y \leq H_i$. Thus Lemma 2.2 holds if we take $r_i(y) =
q_i(1,y)$ and $\alpha_i + \delta$ to be the least value of  $\alpha + M_i\beta$ for which $F_{\alpha,\beta}$ is nonzero other than
$\alpha_i$; each time one takes an $x$ derivative each term in the sum of $(2.6)$ loses a degree in $x$, as does the $O(x^K)$
term. Thus as long as $K$ is chosen sufficiently large (depending on $l$) the conclusions of Lemma 2.2 follow.

\noindent {\heading 3. Proof of Theorem 1.1.}

 Let $k \geq 2$ denote the order of the zero of $S(x,y)$ at the origin. Doing a linear coordinate change if 
necessary, we may assume that for some fixed constant $C_0 > 0$, independent of all other constants in this paper, we have ${1 \over C_0}|\mu_1| < |\mu_2| < C_0 |\mu_1|$ and also that the Taylor expansion $\partial_yS(x,y) = \sum_{\alpha \beta} S_{\alpha \beta}x^{\alpha}y^{\beta}$ 
has a nonvanishing $S_{k-1\,0}x^{k-1}$ term and a nonvanishing $S_{0\,k-1}y^{k-1}$ term, and that the same as true for $\partial_xS(x,y)$. We divide a small rectangle centered at the origin into 8 regions via the lines $y = mx$ and the $x$ and $y$ axes 
as in the beginning of section 2, and then after reflections about the $x$ or $y$ axes and/or the line $y = \pm x$ as necessary we assume we are working on 8 domains of the
form $\{(x,y): 0 < x < b,\, 0 < y < ax\}$. 

We now apply Theorem 2.1 to each $\partial_y S(x,y)$, where $S(x,y)$ now refers to the phase in the possibly reflected coordinates of its domain. Let
$\{D_i\}_{i=1}^n$ denote the domains resulting from applying Theorem 2.1 on these domains; we include the $D_i$ from all 8 domains in a single list. Where $\phi$ is a cutoff function supported on a  small neighborhood of the origin, define $T_i(\lambda_1,\mu_1,\mu_2)$ by
$$T_i(\lambda_1,\mu_1,\mu_2) = \int_{D_i}e^{i\lambda_1S(x,y) +i\mu_1 x + i\mu_2 y}\,\phi(x,y)\,dx\,dy \eqno (3.1)$$
To be perfectly clear, we are still abusing notation slightly in $(3.1)$; $S(x,y)$ denotes the phase function in the reflected coordinates. Since $|\mu_2| \sim |\mu_1|$ (in both the original and reflected coordinates),  to prove Theorem 1.1 it suffices to show that if the support of $\phi$ is sufficiently small, then for each $i$ there is a constant $C$ depending on $S$  such that for $|\mu| > 2$ we have
$$|T_i(\lambda_1,\mu_1,\mu_2)| < C |\mu_2|^{-{1 \over 2}}\|\phi\|_{C^1(V)} \eqno (3.2a)$$
(If $|\mu| \leq 2$ one may just take absolute values and integrate to get the result).
The $i$ for which $\eta_i(x)$ in Theorem 2.1 is of the form $(x, -y + \psi_i(x))$ are dealt with the same way as the $i$
for which $\eta_i(x)$ is of the form $(x, y + \psi_i(x))$, so we always assume $\eta_i(x)$ is of the latter form. 

\noindent Write $T_i(\lambda_1,\mu_1,\mu_2) = T_i^1(\lambda_1,\mu_1,\mu_2) + 
T_i^2(\lambda_1,\mu_1,\mu_2)$, where 
$$T_i^1(\lambda_1,\mu_1,\mu_2) = \int_{\{(x,y) \in D_i:|\mu_2| > 2|\lambda_1 \partial_y S(x,y)|\}} e^{i\lambda_1S(x,y) +i\mu_1 x + i\mu_2 y}\,\phi(x,y)\,dx\,dy \eqno (3.3a)$$
$$T_i^2(\lambda_1,\mu_1,\mu_2) = \int_{\{(x,y) \in D_i:|\mu_2| \leq 2|\lambda_1 \partial_y S(x,y)|\}} e^{i\lambda_1S(x,y) +i\mu_1 x + i\mu_2 y}\,\phi(x,y)\,dx\,dy \eqno (3.3a)$$
We bound $T_i^1(\lambda_1,\mu_1,\mu_2)$ first. We rewrite $(3.3a)$ as

$$T_i^1(\lambda_1,\mu_1,\mu_2) = \int_{\{(x,y) \in D_i:|\mu_2| > 2|\lambda_1 \partial_y S(x,y)|\}} (i\lambda_1\partial_y S(x,y) + i\mu_2)e^{i\lambda_1S(x,y) +i\mu_1 x + i\mu_2 y}$$
$$\times {1 \over i\lambda_1\partial_y S(x,y) + i\mu_2} \,\phi(x,y)\,dx\,dy \eqno (3.4)$$
Note that since $|\mu_2| > 2|\lambda_1 \partial_y S(x,y)|$ in the domain of integration, in the above integration we have 
$|i\lambda_1\partial_y S(x,y) + i\mu_2| > {1 \over 2}|\mu_2|$. This implies that we may integrate by parts in $(3.4)$, integrating  $(i\lambda_1\partial_y S(x,y) + i\mu_2)e^{i\lambda_1S(x,y) +i\mu_1 x + i\mu_2 y}$ and differentiating the other two factors. If the 
derivative lands on $\phi(x,y)$, we take absolute values and integrate, using that $|{1 \over i\lambda_1\partial_y S(x,y) + i\mu_2}|  < 
{2 \over |\mu_2|}$. The result is a bound of $C{1 \over |\mu_2|}\|\phi\|_{C^1(V)}$, a better bound than what we need. If
the derivative lands on ${1 \over i\lambda_1\partial_y S(x,y) + i\mu_2}$, we obtain a term bounded in absolute value by
$$\|\phi\|_{C^1(V)}\int_{\{(x,y) \in D_i:|\mu_2| > 2|\lambda_1 \partial_y S(x,y)|\}}{|\partial_{yy} S(x,y)| \over (\partial_y S(x,y) + \mu_2)^2} \,dx\,dy \eqno (3.5)$$
Because of the linear coordinate change performed at the beginning of the argument, $|\partial_y^k S(x,y)| \neq 0$ on the domain of integration of
$(3.5)$. Thus for fixed $x$, there are boundedly many segments  on which ${|\partial_{yy} S(x,y)| \over (\partial_y S(x,y) + \mu_2)^2} = 
\pm {\partial_{yy} S(x,y) \over (\partial_y S(x,y) + \mu_2)^2}$. On each such segment one can integrate back $\pm {\partial_{yy} S(x,y) \over (\partial_y S(x,y) + \mu_2)^2}$ to obtain $\mp {1 \over \partial_y S(x,y) + \mu_2}$, similar to in the proof of the 
Van der Corput lemma. Since $|{1 \over \partial_y S(x,y) + \mu_2}| \leq 2 {1 \over |\mu_2|}$, we get that $(3.5)$ is bounded
by $C{1 \over |\mu_2|}\|\phi\|_{C^1(V)}$, the same bound as we had for the other term. Lastly, we observe that the endpoint
terms in the integration by parts also give a bound of $C{1 \over |\mu_2|}\|\phi\|_{C^1(V)}$.

We now proceed to bounding $T_i^2(\lambda_1,\mu_1,\mu_2)$. 
The argument from this point on is done somewhat differently if $\beta_i > 0$ or $\beta_i = 0$ for the domain
$D_i$, where $\beta_i$ is as in Theorem 2.1, which we recall we are applying to $\partial_y S(x,y)$.

\noindent {\bf Case 1.} $\beta_i > 0$. We decompose $D_i = \cup_j D_{ijk}$, where $D_{ijk} = \{(x,y) \in D_i: 2^{-j-1} < x \leq 2^{-j}, 2^{-k-1} < y - \psi_i(x) \leq 2^{-k} \}$, and we correspondingly define
$$T_{ijk}^2(\lambda_1,\mu_1,\mu_2) = \int_{\{(x,y) \in D_{ijk}:|\mu_2| \leq 2|\lambda_1 \partial_y S(x,y)|\}}e^{i\lambda_1S(x,y) +i\mu_1 x + i\mu_2 y}\,\phi(x,y)\,dx\,dy \eqno (3.6)$$
The second $y$ derivative of the phase in $(3.6)$ is $\lambda_1 S_{yy}(x,y)$, which by part c) of Theorem 2.1 can be written as
$\lambda_1\beta_id_i x^{\alpha_i} (y - \psi_i(x))^{\beta_i - 1} + o(|\lambda_1 x^{\alpha_i} (y - \psi_i(x))^{\beta_i - 1}|)$. It is here that we use the real-analyticity
condition; if the function is not real-analytic then the error term might not be $o(|\lambda_1 x^{\alpha_i} (y - \psi_i(x))^{\beta_i - 1}|)$ in the event that the lower boundary of $D_i'$ is the $x$-axis. We now apply the measure version of the Van der 
Corput lemma (see [C]) in the $y$ direction, integrate the result in $x$, and we get that
$$|T_{ijk}^2(\lambda_1,\mu_1,\mu_2)| \leq C\|\phi\|_{C^1(V)}2^{-j}(|\lambda_1|^{-{1 \over 2}}2^{{j \alpha_i \over 2}}2^{{k(\beta_i - 1) \over 2}}) \eqno (3.7)$$
In order for $(3.6)$ to be nonzero, there must be at least one point in $D_{ijk}$ for which $|\mu_2| \leq 2|\lambda_1 \partial_y S(x,y)|$. Since $|\lambda_1 \partial_y S(x,y)|$ doesn't vary by more than a constant factor on $D_{ijk}$, this means there 
exists a $C$ such that if $(3.6)$ is nonzero then on all of $D_{ijk}$ one has 
$$ |\mu_2| \leq C|\lambda_1 \partial_y S(x,y)| $$
$$ \leq C'|\lambda_1|2^{-j\alpha - k \beta} \eqno (3.8)$$
Substituting this into $(3.7)$, we get that
 $$|T_{ijk}^2(\lambda_1,\mu_1,\mu_2)| \leq C\|\phi\|_{C^1(V)}2^{-j - {k \over 2}}|\mu_2|^{-{1 \over 2}} \eqno (3.9)$$
We now add $(3.9)$ over all $(j,k)$, resulting in a bound of a constant times $\|\phi\|_{C^1(V)}|\mu_2|^{-{1 \over 2}}$. Since $|\mu_2| \sim |\mu_1|$, this gives us the needed bound of a constant times $\|\phi\|_{C^1(V)}|\mu|^{-{1 \over 2}}$.

\noindent {\bf Case 2.} $\beta_i = 0$. This time we decompose $D_i = \cup_j D_{ij}$ where $D_{ij} = \{(x,y) \in D_i: 2^{-j-1} < x \leq 2^{-j}\}$, and we correspondingly define
$$T_{ij}^2(\lambda_1,\mu_1,\mu_2) = \int_{\{(x,y) \in D_{ij}:|\mu_2| \leq 2|\lambda_1 \partial_y S(x,y)|\}}e^{i\lambda_1S(x,y) +i\mu_1 x + i\mu_2 y}\,\phi(x,y)\,dx\,dy \eqno (3.10)$$
Let $\gamma = {-M_i + \alpha_i \over 2}$.
We write $T_{ij}^2(\lambda_1,\mu_1,\mu_2) = T_{ij}^3(\lambda_1,\mu_1,\mu_2) + T_{ij}^4(\lambda_1,\mu_1,\mu_2)$, where
$$T_{ij}^3(\lambda_1,\mu_1,\mu_2) = \int_{\{(x,y) \in D_{ij}:|\mu_2| \leq 2|\lambda_1 \partial_y S(x,y)|,\,\,|\lambda_1\partial_yS(x,y) + \mu_2| \geq |\lambda_1|^{1 \over 2}2^{-j\gamma} \}}e^{i\lambda_1S(x,y) +i\mu_1 x + i\mu_2 y}$$
$$\times \phi(x,y)\,dx\,dy \eqno (3.11a)$$
$$T_{ij}^4(\lambda_1,\mu_1,\mu_2) = \int_{\{(x,y) \in D_{ij}:|\mu_2| \leq 2|\lambda_1 \partial_y S(x,y)|,\,\, |\lambda_1\partial_yS(x,y) + \mu_2| < |\lambda_1|^{1 \over 2}2^{-j\gamma} \}}e^{i\lambda_1S(x,y) +i\mu_1 x + i\mu_2 y}$$
$$\times \phi(x,y)\,dx\,dy \eqno (3.11b)$$
For $T_{ij}^3(\lambda_1,\mu_1,\mu_2)$ we integrate by parts in $y$ exactly as we did in $(3.4)-(3.5)$, using that $|\lambda_1\partial_yS(x,y) + \mu_2| \geq |\lambda_1|^{1 \over 2}2^{-j\gamma}$ in place of $|\lambda_1\partial_yS(x,y) + \mu_2| \geq {1 \over 2}|\mu_2|$.  Instead of a bound of $C{1 \over |\mu_2|}\|\phi\|_{C^1(V)}$, this time we get the bound
$$ |T_{ij}^3(\lambda_1,\mu_1,\mu_2) | \leq C2^{-j} {1 \over |\lambda_1|^{1 \over 2}2^{-j\gamma}}\|\phi\|_{C^1(V)} \eqno (3.12)$$
$$= C2^{-j}|\lambda_1|^{-{1 \over 2}}2^{j({-M_i + \alpha_i  \over 2})}\|\phi\|_{C^1(V)}\eqno (3.13)$$
Here the ${1 \over |\lambda_1|^{1 \over 2}2^{-j\gamma}}\|\phi\|_{C^1(V)}$ factor is from the $y$ integration and the $2^{-j}$ factor is from the subsequent $x$ integration. 
Like in Case 1, if $(3.11a)$ is nonzero then on the domain of integration we have $|\mu_2| \leq C|\lambda_1 \partial_y S(x,y)|$.
 By Theorem 2.1 c) $|\partial_y S(x,y)| \sim x^{\alpha_i} \sim 2^{-j\alpha_i}$ here (since $\beta_i = 0$). So in $(3.13)$, the $2^{{j \alpha_i \over 2}}$ factor is bounded 
by $C|\lambda_1|^{1 \over 2}|\mu_2|^{-{1 \over 2}}$, and therefore $(3.13)$ is bounded by
$$ C' |\mu_2|^{-{1 \over 2}}2^{j(-{M_i  \over 2} - 1)}\|\phi\|_{C^1(V)}\eqno (3.14)$$
Adding over all $j$ gives a bound of $C'' |\mu_2|^{-{1 \over 2}}\|\phi\|_{C^1(V)}$, the desired bound since $|\mu_2| \sim |\mu_1|$.
We next show that $T_{ij}^4(\lambda_1,\mu_1,\mu_2)$ is also bounded by $(3.13)$, so that $T_{ij}^4(\lambda_1,\mu_1,\mu_2)$ is also bounded by a constant times  $|\mu_2|^{-{1 \over 2}}\|\phi\|_{C^1(V)}$. Taking absolute values in $(3.11b)$ and integrating, we get that $|T_{ij}^4(\lambda_1,\mu_1,\mu_2)|$ is at most
$$\|\phi\|_{C^1(V)} \times  |\{(x,y) \in D_{ij}:|\mu_2| \leq 2|\lambda_1 \partial_y S(x,y)|,\,\, |\lambda_1\partial_yS(x,y) + \mu_2| < |\lambda_1|^{1 \over 2}2^{-j\gamma} \}| \eqno (3.15)$$
We now shift $y$ by $\psi_i(x,y)$, so that where $\eta_i$ is in Theorem 2.1 we have that $|T_{ij}^4(\lambda_1,\mu_1,\mu_2)|$ is at most $\|\phi\|_{C^1(V)}$ times
$$|\{(x,y) \in D_{ij}':|\mu_2| \leq 2|\lambda_1 \partial_y (S \circ \eta_i)(x,y)|,\,\, |\lambda_1\partial_y(S \circ \eta_i)(x,y) + \mu_2| < |\lambda_1|^{1 \over 2}2^{-j\gamma} \}| \eqno (3.16)$$
Here $D_{ij}'$ is the shift of $D_{ij}$ by $\psi_i(x)$ in the $y$ variable. The condition that $|\mu_2| \leq 2|\lambda_1 \partial_y (S \circ \eta_i)(x,y)|$ is used only to go from $(3.13)$ to
$(3.14)$, and we use only the $|\lambda_1\partial_y(S \circ \eta_i)(x,y) + \mu_2| < \lambda_1^{1 \over 2}2^{-j\gamma}$ condition in proving $(3.13)$. So as to be able to use Lemma 2.2, we change
variables from $y$ to $x^{M_i}y$ in $(3.16)$ and get a term bounded by $\|\phi\|_{C^1(V)}$ times
$$2^{-jM_i} |\{(x,y) \in  [2^{-j-1},2^{-j}] \times [0,H_i] : |\lambda_1\partial_y(S \circ \eta_i)(x,x^{M_i}y) + \mu_2| < |\lambda_1|^{1 \over 2}2^{-j\gamma} \}| \eqno (3.17)$$
Our use of $[2^{-j-1},2^{-j}] \times [0,H_i]$ here follows from parts a) and b) of Theorem 2.1. By Lemma 2.2, we have that
$$|\partial_x \big(\lambda_1\partial_y(S \circ \eta_i)(x,x^{M_i}y) + \mu_2)\big)| > C|\lambda_1|x^{\alpha_i - 1}$$
$$ > C'|\lambda_1|2^{-j\alpha_i + j} \eqno (3.18)$$
 Thus for a fixed $y$, the $x$-measure of the set in $(3.17)$ is at most $C|\lambda_1|^{-{1 \over 2}}2^{-j\gamma + j\alpha_i - j}$. Thus $\|\phi\|_{C^1(V)}$ times the quantity in $(3.17)$ is bounded by
$$C|\lambda_1|^{-{1 \over 2}}2^{-j\gamma + j\alpha_i - j -jM_i} \|\phi\|_{C^1(V)}\eqno (3.19)$$
Substituting back in for $\gamma$, this becomes
 $$C|\lambda_1|^{-{1 \over 2}}2^{-j{M_i \over 2} + j{\alpha_i \over 2}- j} \|\phi\|_{C^1(V)}\eqno (3.20)$$
This is exactly $(3.13)$. The condition that $|\mu_2| \leq 2|\lambda_1 \partial_y S(x,y)|$ is now 
used exactly as it was when going from $(3.13)$ to $(3.14)$. This again leads to the bound $(3.14)$ for $|T_{ij}^4(\lambda_1,\mu_1,\mu_2)|$, and after summing this in $j$ we are done.

\line{}\line{}

\noindent {\heading 4. Proof of Theorem 1.2.}

In the proof of Theorem 1.2 we will make use of sublevel set estimates that are analogous to the oscillatory integral estimates we have
been using. Specifically, if $f(x,y)$ is real analytic on a neighborhood of the origin such that $f(0,0) = 0$ and $\nabla f(0,0) = 0$, 
for a given $U$ contained in the domain of $f(x,y)$ and an $0 < r < {1 \over 2}$ we define
$$A_U(r) = |\{(x,y) \in U: |f(x,y)| < r \}| \eqno (4.1)$$
Using resolution of singularities (see [AGV] Ch. 6 for details), in the real-analytic case if $U$ is a sufficiently small ball centered 
at the origin then  as $r \rightarrow 0$ one has an asymptotic expansion of the form
$$A_U(r) = C_U r^{\epsilon} |\ln(r)|^m + o( r^{\epsilon} |\ln(r)|^m) \eqno (4.2)$$
Here $C_U > 0$ and $(\epsilon,m)$ is the same as in $(1.2)$, unless $(\epsilon,m) = (1,0)$, in which case $(\epsilon,m)$ could be $(1,0)$ or $(1,1)$. In [G1] it is shown that in the general smooth case,
an analogue of $(4.2)$ holds. Namely, there is a $C_U$ such that $A_U(r) \leq  C_U r^{\epsilon} |\ln(r)|^m$, and often $(4.2)$ still holds. In [G1] it is shown that in the cases where $(4.2)$ does not hold, $m$ is always $0$ and for all $\epsilon ' > \epsilon$ there is a 
constant $C_{U,\epsilon'} > 0$ such that $A_U(r) \geq  C_{U,\epsilon'} r^{\epsilon'}$. This extension to the smooth case does
use the notion of adapted coordinate systems, and is the only way in which this paper relies on them. However, one can avoid 
relying on the 
use of adapted coordinate systems entirely by doing arguments very similar to those of [G1] on the constructions of Theorem 2.1.

\noindent The above discussion leads to the following lemma.

\noindent {\bf Lemma 4.1.} Let $(\epsilon,m)$ be as above, and let $\{D_i'\}_{i=1}^n$ be the domains obtained by applying Theorem 2.1 to $f(x,y)$, and let $(\alpha_i,\beta_i)$ be as in that theorem. Then there exists a constant $C$ such that for each $i$ and all $0 < r < {1 \over 2}$ we have
$$|\{(x,y) \in D_i': x^{\alpha_i}y^{\beta_i} < r\}| \leq C r^{\epsilon}|\ln(r)|^m \eqno (4.3)$$

\noindent {\bf Proof.} In the case that $g_i(x)$ is not identically zero in Theorem 2.1, a sliver $\{(x,y): 0 < x < b: 0 < y < Cx^{M_i}\}$ is disjoint from $D_i'$, so $x^{\alpha_i}y^{\beta_i}$ is bounded below by $C'x^{\alpha_i - M_i\beta_i}$ on $D_i'$. In part c) of Theorem 2.1, since $f_2$ has a zero of infinite order at the origin, for any $N$ one has an 
estimate of the form $|f_2(x,y)| < C_N x^N$. Thus in $(2.2)$ one can replace $f_1$ by $f$, which implies $f \circ \eta_i$ is
within a constant factor of $x^{\alpha_i}y^{\beta_i}$ on $D_i$. Since the Jacobian of $\eta_i$ is everywhere equal to $1$, the
measure of the sublevel sets of $|f \circ \eta_i|$ will be no greater than the measure of the corresponding sublevel sets of $|f|$. Thus $(4.3)$ holds.

Now consider the case where $g_i(x)$ is identically zero. Define $D_i^N = \{(x,y) \in D_i': y > x^N\}$. Then exactly as above
one has that $|\{(x,y) \in D_i^N: x^{\alpha_i}y^{\beta_i} < r\}| \leq C_N r^{\epsilon}|\ln(r)|^m$. In the case that $\alpha_i \geq \beta_i$
this is enough; a direct calculation reveals that for large enough $N$, $|\{(x,y) \in D_i^N: x^{\alpha_i}y^{\beta_i} < r\}| $ is 
within a constant factor of $|\{(x,y) \in D_i': x^{\alpha_i}y^{\beta_i} < r\}|$. In the case where $\alpha_i < \beta_i$, a direct
calculation reveals 
that $|\{(x,y) \in D_i^N: x^{\alpha_i}y^{\beta_i} < r\}|$ is of the form $C_N r^{\delta_N} +$ lower order terms, and that
$|\{(x,y) \in D_i': x^{\alpha_i}y^{\beta_i} < r\}|$ is of the form $Cr^{\delta} +$ lower order terms where $\lim_{N \rightarrow
\infty} \delta_N = \delta$. So like in the previous paragraph, $\delta_N \geq \epsilon$ for each $N$, and taking limits as $N \rightarrow \infty $ we get that $\delta \geq \epsilon$. So regardless of whether or not $m=0$ or $1$, $(4.3)$ will hold and we are done.

\noindent We now are in a position to prove Theorem 1.2.

\noindent {\bf Proof of Theorem 1.2}. 

Let $k > 0$ be the order of the zero of $S(x,y)$ at the origin. Rotating coordinates if
necessary, we assume that the Taylor expansion $\sum_{\alpha \beta} S_{\alpha \beta}x^{\alpha}y^{\beta}$ of $S(x,y)$ has nonvanishing $S_{0 k}y^k$ and $S_{k 0} x^k$ terms. We perform the reflections at the beginning of section 2 and then apply Theorem 2.1 to (the reflected) $S(x,y)$.
Note this is a different function from the previous section. Let $\{D_i\}_{i=1}^n$ be all of the resulting regions. We will bound the
portion of $T(\lambda_1,\lambda_2,\lambda_3) = T(\lambda_1,\mu_1,\mu_2)$ coming from a given $D_i$ and sum over all $i$.
We will slightly abuse notation in the following and refer to a reflected $S(x,y)$ as just $S(x,y)$.
The argument
naturally breaks into three cases. The first is when $\beta_i$  in Theorem 2.1 is  greater than $1$ and the lower boundary of $D_i'$ is the $x$-axis (in other words, $g_i(x)$ is identically zero), the second is when either $\beta_i >  1$ and the lower boundary of $D_i'$ is not the $x$-axis, or $\beta_i = 1$ (with no restrictions), and the third case is when $\beta_i$ is zero.

\noindent {\bf Case 1.} $\beta_i > 1$ and $g_i(x)$ is identically zero. 

Consider $(2.2)$ when $l = 0$ and $m = \beta_i$. Because
in $(2.2)$ the function $f_2$  has a zero of infinite order at the origin, on $D_i'$ one has 
$|\partial_y^{\beta_i}(f_2 \circ \eta_i)(x,y)| < C_N x^N$ for any $N$. Thus we may replace $f_1$ by $f$ (which is $S$ here)
to obtain  that for some constant $C$  we have 
$$ |\partial_y^{\beta_i}(S \circ \eta_i)(x,y)| > C x^{\alpha_i} \eqno (4.4)$$
Denote by $T_i(\lambda_1,\mu_1,\mu_2)$  portion of the integral $(1.1)$ coming from $D_i$. In this integral, we do the coordinate
change $\eta_i(x,y)$ given by Theorem 2.1, obtaining
$$ T_i(\lambda_1,\mu_1,\mu_2) = \int_{D_i'} e^{i\lambda_1 (S\circ \eta_i)(x,y) + i\mu_1 x +  i\mu_2y + i\mu_2\psi_i(x)} \phi_i(x,y)
\,dx\,dy \eqno (4.5)$$
The $ i\mu_2 y$ term might have a minus sign in front, but since that case is done
exactly the same way we will assume $T_i(\lambda_1,\mu_1,\mu_2)$ is of the form $(4.5)$. Here $\phi_i(x,y)$ is a compactly supported function such that $\phi_i(x^N,y)$ is smooth for some $N$.
We now dyadically decompose $T_i = \cup_j T_{ij}$. Denoting by $D_{ij}$ the set $\{(x,y) \in D_i': 2^{-j-1} \leq x < 2^{-j}\}$, we define
$$ T_{ij}(\lambda_1,\mu_1,\mu_2) = \int_{D_{ij}} e^{i\lambda_1 (S\circ \eta_i)(x,y) + i\mu_1 x +  i\mu_2y + i\mu_2\psi_i(x)}  \phi_i(x,y)\,dx\,dy \eqno (4.6)$$
We now apply the standard Van der Corput lemma (see [S] ch 8) in $(4.6)$ in the $y$ direction, using $(4.4)$, and then integrate the result in $x$.
One gets
$$ |T_{ij}(\lambda_1,\mu_1,\mu_2)| \leq C \|\phi\|_{C^1(V)}|\lambda_1|^{-{1 \over \beta_i}}2^{j {\alpha_i \over \beta_i}}  \times 
2^{-j} \eqno (4.7)$$
In $(4.6)$ we can get a crude estimate by taking absolute values and integrating, obtaining a constant times $\|\phi\|_{C^1(V)}2^{-j-jM_i}$, $M_i$ as in Theorem 2.1a). Thus one can extend $(4.7)$ to
$$ |T_{ij}(\lambda_1,\mu_1,\mu_2)| \leq C\|\phi\|_{C^1(V)}\min(2^{-j-jM_i},  |\lambda_1|^{-{1 \over \beta_i}}2^{j {\alpha_i \over \beta_i}}  \times 
2^{-j}) \eqno (4.8)$$
An elementary calculation reveals that the measure of $\{(x,y) \in D_{ij}: x^{\alpha_i}y^{\beta_i} < {1 \over |\lambda_1|}\}$ is within a 
constant factor of $\min(2^{-j-jM_i},  |\lambda_1|^{-{1 \over \beta_i}}2^{j {\alpha_i \over \beta_i}}  \times 
2^{-j})$. Thus we have
$$|T_{ij}(\lambda_1,\mu_1,\mu_2)| \leq C'\|\phi\|_{C^1(V)}\big|\{(x,y) \in D_{ij}: x^{\alpha_i}y^{\beta_i} < |\lambda_1|^{-1}\}\big| \eqno (4.9)$$
Adding over all $j$ then gives
$$|T_i(\lambda_1,\mu_1,\mu_2)| \leq C'\|\phi\|_{C^1(V)}\big|\{(x,y) \in D_i': x^{\alpha_i}y^{\beta_i} <|\lambda_1|^{-1}\}\big| \eqno (4.10)$$
By Lemma 4.1, the right-hand side of $(4.10)$ is bounded by $C\|\phi\|_{C^1(V)}|\lambda_1|^{-\epsilon} (\ln |\lambda_1|)^m$ and we are done. Note that the Case 1 argument did not use any restrictions on the value of $\epsilon$.

\noindent {\bf Case 2.} $\beta_i = 1$ or $\beta_i > 1$ and $g_i(x)$ is not identically zero. 

 Similar to in the previous case, when $l = m = 1$ we can replace $f_1$ by $f$ in $(2.2)$. When $\beta_i = 1$ this is because $x^{\alpha_i - 1}y^{\beta_i - 1}$ is a power of $x$, so since
$|{\partial^2 \over \partial x \partial y}(f_2 \circ \eta_i)(x,y)| < C_N x^N$ for any given $N$, changing from $f_1$ to $f$ (which is $S$ here) will
not interfere with the validity of $(2.2)$. When $\beta_i > 1$ and $g_i(x)$ is not identically zero, we may do this replacement since
a sliver $\{(x,y): 0 < x < b: 0 < y < Cx^{M_i}\}$ is disjoint from the domain, so that $x^{\alpha_i - 1}y^{\beta_i -1}$ is 
bounded below by $C'x^{\alpha_i - 1 - M_i\beta_i - M_i}$ and similar considerations apply. Thus we may replace $f_1$ by $f = S$ in
$(2.2)$, and on $D_i'$  we have a lower bound of the form
$$\bigg|{\partial^2 \over \partial x \partial y}(S\circ \eta_i)(x,y)\bigg| > C x^{\alpha_i - 1}y^{\beta_i - 1} \eqno (4.11)$$
Note that ${\partial^2 \over \partial x \partial y} (i\lambda_1 (S\circ \eta_i)(x,y) + i\mu_1 x +  i\mu_2y + i\mu_2\psi_i(x))  = 
 i\lambda_1{\partial^2 \over \partial x \partial y} (S\circ \eta_i)(x,y)$. Thus $(4.11)$ is relevant to  $T_{ijk}(\lambda_1,\mu_1,\mu_2)$. 

This time we dyadically decompose $(4.5)$ in both the $x$ and $y$ directions. Namely, let $D_{ijk} = \{(x,y) \in D_i': 2^{-j-1} \leq x < 2^{-j},\,\,2^{-k-1} \leq y < 2^{-k} \}$ and define $T_{ijk}(\lambda_1,\mu_1,\mu_2)$ by 
$$ T_{ijk}(\lambda_1,\mu_1,\mu_2) = \int_{D_{ijk}} e^{i\lambda_1 (S\circ \eta_i)(x,y) + i\mu_1 x +  i\mu_2y + i\mu_2\psi_i(x)}  \phi_i(x,y)\,dx\,dy \eqno (4.12)$$

We now proceed similarly to in Case 2 of Theorem 1.1. We write  $T_{ijk}(\lambda_1,\mu_1,\mu_2)  =  T_{ijk}^1(\lambda_1,\mu_1,\mu_2) +  T_{ijk}^2(\lambda_1,\mu_1,\mu_2)$, where 
$$ T_{ijk}^1(\lambda_1,\mu_1,\mu_2) = \int_{\{(x,y) \in D_{ijk}: |\lambda_1\partial_yS + \mu_2| < |\lambda_1|^{{1 \over 2}}2^{k -{j\alpha_i + k\beta_i \over 2}}\}} e^{i\lambda_1 (S\circ \eta_i)(x,y) + i\mu_1 x +  i\mu_2y + i\mu_2\psi_i(x)}$$
$$\times \phi_i(x,y)\,dx\,dy \eqno (4.13a)$$
$$ T_{ijk}^2(\lambda_1,\mu_1,\mu_2) = \int_{\{(x,y) \in D_{ijk}: |\lambda_1\partial_yS + \mu_2| \geq  |\lambda_1|^{{1 \over 2}}2^{k -{j\alpha_i + k\beta_i \over 2}}\}} e^{i\lambda_1 (S\circ \eta_i)(x,y) + i\mu_1 x +  i\mu_2y + i\mu_2\psi_i(x)}$$
$$\times \phi_i(x,y)\,dx\,dy \eqno (4.13b)$$
For $(4.13a)$, we simply take absolute values and integrate, obtaining
$$ |T_{ijk}^1(\lambda_1,\mu_1,\mu_2)| \leq C|\{(x,y) \in D_{ijk}: |\lambda_1\partial_yS + \mu_2| <  |\lambda_1|^{{1 \over 2}}2^{k -{j\alpha_i + k\beta_i \over 2}}\}| \,\,\|\phi\|_{C^1(V)}
\eqno (4.14)$$
By $(4.11)$, the absolute value of the $x$ derivative of $\lambda_1\partial_yS + \mu_2$ is bounded below by $C|\lambda_1|2^{-j(\alpha_i - 1) - k(\beta_i - 1)}$. So for a 
given $y$, the $x$-measure of the set in $(4.14)$ is at most $C|\lambda_1|^{-{1 \over 2}}2^{k -{j\alpha_i + k\beta_i \over 2}}\times 2^{j(\alpha_i - 1) + k(\beta_i - 1)} = C|\lambda_1|^{-{1 \over 2}}2^{-j + {j \alpha_i+ k\beta_i\over 2}}$. Inserting this into $(4.14)$ and then integrating in $y$, we get that
$$|T_{ijk}^1(\lambda_1,\mu_1,\mu_2)| \leq C |\lambda_1|^{-{1 \over 2}}2^{-j - k  + {j\alpha_i + k\beta_i \over 2}} \|\phi\|_{C^1(V)}\eqno (4.15)$$
This is the estimate we will need. Moving on to $ T_{ijk}^2(\lambda_1,\mu_1,\mu_2) $, we integrate by parts in $y$ in $(4.13b)$.
We write $e^{i\lambda_1 (S\circ \eta_i)(x,y) + i\mu_1 x +  i\mu_2y + i\mu_2\psi_i(x)}$ as
$$( {i\lambda_1\partial_y(S \circ \eta_i) (x,y)+ i\mu_2}) e^{i\lambda_1 (S\circ \eta_i)(x,y) + i\mu_1 x +  i\mu_2y + i\mu_2\psi_i(x)}\times {1 \over {i\lambda_1\partial_y(S \circ \eta_i)(x,y)+ i\mu_2}} \eqno (4.16)$$
We integrate by parts in $(4.13b)$
by integrating the left factor of $(4.16)$ and differentiating the rest. There are two places the derivative may land, the ${1 \over {i\lambda_1\partial_y(S \circ \eta_i)(x,y)+ i\mu_2}}$ factor and $\phi_i(x,y)$. In the first case, the differentiation gives ${-\lambda_1\partial_{yy}(S \circ \eta_i)(x,y) \over i(\partial_y(S \circ \eta_i)(x,y)+ \mu_2)^2}$. We then take absolute values and integrate in $y$, very similar to in $(3.5)$. As in that situation, the end result of
the $y$ integration is a bound of $C\|\phi\|_{C^1(V)}$ times the maximum of $|{1 \over {i\lambda_1\partial_y(S \circ \eta_i)(x,y)+ i\mu_2}}|$ on the
domain of integration, or $C |\lambda_1|^{-{1 \over 2}}2^{-k +{j\alpha_i + k\beta_i \over 2}}\|\phi\|_{C^1(V)}$. We then do
the $x$ integration to get an overall factor of $C |\lambda_1|^{-{1 \over 2}}2^{-j - k +{j\alpha_i + k\beta_i \over 2}}\|\phi\|_{C^1(V)}$, which is the same as in $(4.15)$.

The second place the derivative may land is the $\phi_i(x,y)$ term. In this case we take absolute values, bound $|{1 \over {i\lambda_1\partial_yS(x,y)+ i\mu_2}}|$ by  $|\lambda_1|^{-{1 \over 2}}2^{-k +{j\alpha_i + k\beta_i \over 2}}$, and integrate.
The result is $C |\lambda_1|^{-{1 \over 2}}2^{-j - 2k +{j\alpha_i + k\beta_i \over 2}}\|\phi\|_{C^1(V)}$, better than what we 
need.  Lastly, we have the endpoint terms of the integration by parts, which will give the same bounds as in the last paragraph, or $C |\lambda_1|^{-{1 \over 2}}2^{-j - k +{j\alpha_i + k\beta_i \over 2}}\|\phi\|_{C^1(V)}$.

\noindent Putting the above together, we conclude that
$$|T_{ijk}(\lambda_1,\mu_1,\mu_2)| \leq C |\lambda_1|^{-{1 \over 2}}2^{-j - k +{j\alpha_i + k\beta_i \over 2}}\|\phi\|_{C^1(V)} \eqno (4.17)$$
We can rewrite this as
$$|T_{ijk}(\lambda_1,\mu_1,\mu_2)| \leq C' \|\phi\|_{C^1(V)} \int_{D_{ijk}}{1 \over (|\lambda_1|x^{\alpha_i}y^{\beta_i})^{1 \over 2}} \eqno (4.18)$$
One can take absolute values in $(4.12)$ and integrate to get another (crude) bound for $|T_{ijk}(\lambda_1,\mu_1,\mu_2)|$. Incorporating this into $(4.18)$ gives
$$|T_{ijk}(\lambda_1,\mu_1,\mu_2)| \leq C' \|\phi\|_{C^1(V)} \int_{D_{ijk}}\min\bigg(1, {1 \over (|\lambda_1|x^{\alpha_i}y^{\beta_i})^{1 \over 2}}\bigg) \eqno (4.19)$$
Adding this over all $j$ and $k$ then gives
$$|T_i(\lambda_1,\mu_1,\mu_2)| \leq C' \|\phi\|_{C^1(V)} \int_{D_i'}\min\bigg(1, {1 \over (|\lambda_1|x^{\alpha_i}y^{\beta_i})^{1 \over 2}}\bigg) \eqno (4.20)$$
$$ = C' \|\phi\|_{C^1(V)}\bigg(\big|\{(x,y) \in D_i': x^{\alpha_i}y^{\beta_i} < {1 \over |\lambda_1|}\}\big| + {1 \over |\lambda_1|^{1 \over 2}}
\int_{\{(x,y) \in D_i':\,\, x^{\alpha_i}y^{\beta_i} \geq {1 \over |\lambda_1|}\}} {1 \over (x^{\alpha_i}y^{\beta_i})^{1 \over 2}}\bigg)\eqno (4.21)$$
By Lemma 4.1, the measure of $\{(x,y) \in D_i': x^{\alpha_i}y^{\beta_i} < {1 \over |\lambda_1|}\}$ is bounded by 
$C|\lambda_1|^{-\epsilon}\ln|\lambda_1|$, the desired estimate. As for the second term, we write the integral in terms 
of distribution functions. Namely, we have
$$\int_{\{(x,y) \in D_i': x^{\alpha_i}y^{\beta_i} \geq {1 \over |\lambda_1|}\}} {1 \over (x^{\alpha_i}y^{\beta_i})^{1 \over 2}}
= \int_{1 \over |\lambda_1|}^{\infty} {1 \over 2} t^{-{3 \over 2}}|\{(x,y) \in D_i': {1 \over |\lambda_1|} \leq x^{\alpha_i}y^{\beta_i} \leq t\}|\,dt\eqno (4.22)$$
By Lemma 4.1, this is at most
$$ C\int_{1 \over |\lambda_1|}^1 {1 \over 2} t^{-{3 \over 2}}t^{\epsilon}\ln(t)^m \,dt \eqno (4.23)$$
Since $\epsilon$ is being assumed to be at most ${1 \over 3}$ here, the integral converges and is bounded by $C|\lambda_1|^{{1 \over 2} - \epsilon}(\ln |\lambda_1|)^m$. Thus the second term of $(4.21)$ is bounded by 
$$C\|\phi\|_{C^1(V)}|\lambda_1|^{-{1 \over 2}}|\lambda_1|^{{1 \over 2} - \epsilon}(\ln |\lambda_1|)^m$$
$$= C\|\phi\|_{C^1(V)}|\lambda_1|^{-\epsilon} (\ln |\lambda_1|)^m \eqno (4.24)$$
This is the desired estimate and we are done. Note that the Case 2 argument only required that $\epsilon < {1 \over 2}$, which was
used to say that $(4.23)$ converges.

\noindent {\bf Case 3.} $\beta_i = 0$.

For this case, it will be helpful to use the following consequence of the Van der Corput lemma 
from [ArCuKa] which was also used in [IKeM][IM2].

\noindent {\bf Lemma 4.2.} Suppose $f$ is a smooth real-valued function on an interval $I$ such that for some integer $n \geq 2$ and some constants $C, C' > 0$, for all $t \in I$ one has that $C' \leq \sum_{i = 2}^n |f^{(i)}(t)| \leq C$. Then there is a constant
$C''$ depending only on $C$ and $C'$ such that for all $\lambda \in \R$ one has
$$\bigg|\int_I e^{i\lambda f(t)} \phi(t)\,dt\bigg| \leq C''(\|\phi\|_{L^{\infty}(I)} + \|\phi'\|_{L^1(I)})(1 + |\lambda|)^{-{1 \over n}} \eqno (4.25)$$
We now start the Case 3 argument. By theorem 2.1a)-b), when $\beta_i = 0$ the domain $D_i'$ is of the form $\{(x,y): 0 < x < b,\,\,0 < y < H_ix^{M_i}\}$. Thus we have
$$ T_i(\lambda_1,\mu_1,\mu_2) = \int_{\{(x,y): 0 < x < b,\,\,0 < y < H_ix^{M_i}\}}  e^{i\lambda_1 (S\circ \eta_i)(x,y) + i\mu_1 x +  i\mu_2y + i\mu_2\psi_i(x)}\phi_i(x,y) \,dx\,dy \eqno (4.26)$$
By part d) of Theorem 2.1, either $\psi_i(x) = k_ix$ for some $k_i$ (possibly zero), or $\psi_i(x) = k_ix + l_ix^{s_i} +$
higher order terms (if any), where $l_i \neq 0$ and $1 < s_i \leq M_i$. Thus if we write $\xi(x) = \psi_i(x) - k_ix$, then
either $\xi(x) = 0$ or $\xi(x)$ as a zero of order $s_i \leq M_i$ at the origin. Letting $\mu_3 = \mu_1 + k_i\mu_2$, we 
correspondingly write the expression $(4.26)$ for  $T_i(\lambda_1,\mu_1,\mu_2)$ as
$$ \int_{\{(x,y): 0 < x < b,\,\,0 < y < H_ix^{M_i}\}}  e^{i\lambda_1 (S\circ \eta_i)(x,y) + i\mu_3 x + i\mu_2(\xi(x) + y)}\phi_i(x,y) \,dx\,dy \eqno (4.27)$$

By Theorem 2.1b) there is some $(\alpha,\beta) \neq (\alpha_i,0)$,  $\beta$ a positive integer,  for which the Taylor expansion of $S \circ \eta_i (x,y)$ has a nonzero $S_{\alpha\beta} x^{\alpha}y^{\beta}$ term and for which $\alpha + M_i\beta = \alpha_i$. So $\alpha_i \geq M_i$. Since
$S(x,y)$ is assumed to have a zero of order at least 2 at the origin, 
$(\alpha, \beta)$ cannot be $(0,1)$ and the statement $\alpha + M_i\beta = \alpha_i$ in fact implies that $\alpha_i >  M_i$.

We first prove Theorem 1.2 under the assumption that $\xi(x)$ is not identically zero and $\epsilon < {1 \over 3}$; the modifications needed in the cases where $\epsilon = {1 \over 3}$ or where $\xi(x)$ is identically zero will be described afterwards.

So assuming $\xi(x)$ is not identically zero, for fixed $y$ the phase in $(4.27)$ can be written as the sum of three terms. The first is $\lambda_1 (S\circ \eta_i)(x,y)$, which is of the form $\lambda_1 d_i x^{\alpha_i}$ plus a small error term by Theorem 2.1b), with corresponding expressions for its $x$ derivatives. The second term is $\mu_2(\xi(x) + y)$, where $\mu_2 \xi(x)$ is of the form $\mu_2\big(l_ix^{s_i} + O(x^{s_i + \delta})\big)$, and the third is $\mu_3x$. It is the first two terms that concern us here. Note that $\alpha_i > M_i  \geq s_i > 1$, so that the exponents $\alpha_i$ and $s_i$ are distinct.  As a result, the 2 by 2 matrix $A_i$
with rows $(\alpha_i(\alpha_i - 1), s_i(s_i - 1))$ and $(\alpha_i(\alpha_i - 1)(\alpha_i - 2), s_i(s_i - 1)(s_i - 2))$ has determinant $\alpha_i(\alpha_i -1)s_i(s_i - 1)(\alpha_i - s_i) \neq 0$. Thus for some constants $c$ and $c'$, for any vector $v$ one has $  c'\|v\| \geq \|A_iv\| \geq c\|v\|$.
In particular, letting $v =  ( \lambda_1d_ix^{\alpha_i}, \mu_2l_ix^{s_i})$, we have 
$$ 
|\lambda_1\alpha_i(\alpha_i - 1)d_ix^{\alpha_i} + s_i(s_i-1)\mu_2l_ix^{s_i}| +  |\lambda_1\alpha_i(\alpha_i - 1)(\alpha_i - 2)d_ix^{\alpha_i} + s_i(s_i-1)(s_i-2)\mu_2l_ix^{s_i}|$$
$$\geq c|(\lambda_1 d_ix^{\alpha_i}, \mu_2l_ix^{s_i})|
\eqno (4.28)$$
 Restating, $(4.28)$ implies that for $f(x) = d_ix^{\alpha_i} +  \mu_2l_ix^{s_i} +  \mu_3x$ we have
$$ |x^2f''(x)| + |x^3f'''(x)| \geq c'( |\lambda_1d_ix^{\alpha_i}| +  |\mu_2l_ix^{s_i}| )\eqno (4.29)$$
Adjusting for error terms, if $x$ is sufficiently small and the $\delta$ coming from Theorem 2.1b)  is sufficiently small (recall it can be chosen independent of the $s_i$ and $\alpha_i$), then  independent of the parameters $\lambda_1, \mu_2,$ and $\mu_3$, if $p_y(x)$ denotes the phase function $\lambda_1 (S\circ \eta_i)(x,y) + \mu_3 x + \mu_2(\xi(x) + y)$
we similarly have
$$ |x^2p_y''(x)| + |x^3p_y'''(x)| \geq c''(|\lambda_1d_ix^{\alpha_i}| +  |\mu_2l_ix^{s_i}|)\eqno (4.30)$$
Next we dyadically decompose $(4.27)$, writing $T_i = \sum_j T_{ij}$, where $T_{ij}$ is defined by
$$ T_{ij}(\lambda_1,\mu_1,\mu_2) =  \int_{\{(x,y): 2^{-j-1}  < x < 2^{-j},\,\,0 < y < H_ix^{M_i}\}}  e^{ip_y(x)} \phi_i(x,y) \,dx\,dy \eqno (4.31)$$
We scale $(4.31)$ in $x$, obtaining
$$ T_{ij}(\lambda_1,\mu_1,\mu_2) =  2^{-j}\int_{\{(x,y): {1 \over 2}  < x < 1,\,\,0 < y < H_ix^{M_i}\}}  e^{ip_y(2^{-j}x)} \phi_i(2^{-j}x,y) \,dx\,dy \eqno (4.32)$$
By $(4.30)$, the phase function $q(x) = p_y(2^{-j}x)$ satisfies
$$ |q''(x)| + |q'''(x)| \geq c( |\lambda_12^{-j\alpha_i}| + |\mu_2 2^{-js_i}|) \eqno (4.33)$$
We now apply the Van der Corput-type lemma, Lemma 4.2, in the $x$ direction letting $f(t) =  (|\lambda_12^{-j\alpha_i}| + |\mu_2 2^{-js_i}| )q(x)$, and letting $|\lambda_12^{-j\alpha_i}| + |\mu_2 2^{-js_i}|$ be what is called $\lambda$ in that lemma. The cutoff function
$\phi_i(2^{-j}x,y)$ is equal to $\phi(2^{-j}x,y + \psi_i(2^{-j}x))$ where $\psi_i(x)$ is of
the form $\zeta(x^{1 \over N})$ for a smooth $\zeta$, for some large $N$. Thus the effect of this cutoff function in an application
of Lemma 4.2 in the $x$ direction is to an introduce a factor bounded by (something slightly better than) $C\|\phi\|_{C^1(V)}$. After applying Lemma 4.2 in the
$x$ direction with $n = 3$ and then integrating in $y$ we get
$$|T_{ij}(\lambda_1,\mu_1,\mu_2)| \leq C2^{-j(M_i+1)}\|\phi\|_{C^1(V)}(|\lambda_12^{-j\alpha_i}| + |\mu_2 2^{-js_i}|)^{-{1 \over 3}}\eqno (4.34)$$
We are only interested in the first of the two  terms in the right hand side of $(4.34)$, so we use the bound 
$$|T_{ij}(\lambda_1,\mu_1,\mu_2)| \leq C\|\phi\|_{C^1(V)}2^{-j(M_i+1)}|\lambda_12^{-j \alpha_i}|^{-{1 \over 3}} \eqno (4.35)$$
By simply taking absolute values and integrating in $(4.31)$, one has
$$|T_{ij}(\lambda_1,\mu_1,\mu_2)| \leq C2^{-j(M_i+1)}\|\phi\|_{C^1(V)}$$
Combining this with $(4.35)$, we get that
$$|T_{ij}(\lambda_1,\mu_1,\mu_2)| \leq C2^{-j(M_i+1)}\|\phi\|_{C^1(V)}\min(1, |\lambda_12^{-j \alpha_i}|^{-{1 \over 3}}) \eqno (4.36)$$
By Lemma 2.2, $|S \circ \eta(x,y)| \sim x^{\alpha_i} \sim 2^{-j\alpha_i}$ on $D_i'$, and furthermore the portion of $D_i'$ between $2^{-j-1}$ and $2^{-j}$
has measure $\sim 2^{-j} \times 2^{-jM_i}$. So $(4.36)$ implies that 
$$|T_{ij}(\lambda_1,\mu_1,\mu_2)| \leq C'\|\phi\|_{C^1(V)}\int_{\{(x,y) \in D_i': \,\,x \in [2^{-j-1}, 2^{-j}]\}}\min(1, |\lambda_1x^{\alpha_i}|^{-{1 \over 3}}) \eqno (4.37)$$
We now argue as in $(4.19)-(4.24)$. Adding $(4.37)$ over all $j$ gives
$$|T_i(\lambda_1,\mu_1,\mu_2)| \leq C'\|\phi\|_{C^1(V)}\int_{D_i'}\min(1, |\lambda_1x^{\alpha_i}|^{-{1 \over 3}}) \eqno (4.38)$$
$$ = C' \|\phi\|_{C^1(V)}\bigg(\big|\{(x,y) \in D_i': x^{\alpha_i} < {1 \over |\lambda_1|}\}\big| + {1 \over |\lambda_1|^{1 \over 3}}
\int_{\{(x,y) \in D_i': x^{\alpha_i} \geq {1 \over |\lambda_1|}\}} {1 \over (x^{\alpha_i})^{1 \over 3}}\bigg)\eqno (4.39)$$
For the first term in $(4.39)$, by Lemma 4.1 the measure of $\{(x,y) \in D_i': x^{\alpha_i} < {1 \over |\lambda_1|}\}$ is bounded by 
$C|\lambda_1|^{-\epsilon}\ln|\lambda_1|$, the desired estimate. In view of the form of this sublevel set, we won't have
a logarithmic factor so we even have that 
$$|\{(x,y) \in D_i': x^{\alpha_i} < {1 \over |\lambda_1|}\}| \leq C|\lambda_1|^{-\epsilon} \eqno (4.40)$$
For the second term of $(4.39)$, like before we write the integral in terms 
of distribution functions. We get
$$\int_{\{(x,y) \in D_i': x^{\alpha_i} \geq {1 \over |\lambda_1|}\}} {1 \over (x^{\alpha_i})^{1 \over 3}}
= \int_{1 \over |\lambda_1|}^{\infty} {1 \over 3} t^{-{4 \over 3}}|\{(x,y) \in D_i': {1 \over |\lambda_1|} \leq x^{\alpha_i} \leq t\}|\,dt$$
Inserting $(4.40)$, this is at most
$$ C\int_{1 \over |\lambda_1|}^{\infty} {1 \over 3} t^{-{4 \over 3}}t^{\epsilon} \,dt \eqno (4.41)$$
Since we are assuming $\epsilon < {1 \over 3}$ for now, the integral $(4.41)$ is absolutely integrable and is bounded by $C'|\lambda_1|^{-\epsilon + {1 \over 3}}$. Thus the second term in the parentheses of $(4.39)$ is bounded by 
 $C'|\lambda_1|^{-\epsilon}$. Adding together with the first term we see that we have the
 estimate
$$|T_i(\lambda_1,\mu_1,\mu_2)| \leq C \|\phi\|_{C^1(V)}|\lambda_1|^{-\epsilon} \eqno (4.42)$$
This gives the desired bounds (and even an extra logarithm when $m$ = 1).
This concludes the proof where $\epsilon < {1 \over 3}$ and $\xi(x)$ does not have a zero of infinite order at $x = 0$.

We now assume $\xi(x)$ is identically zero, for any $\epsilon \leq {1 \over 3}$. Then for fixed $y$ the phase in $(4.27)$ is the sum
$\lambda_1 (S\circ \eta_i (x,y))$ and a linear function of $x$. Thus by $(2.1)$ we can just use Lemma 4.2 for 
second derivatives here, since the second $x$-derivative of a linear function is zero.
So if $p_y(x)$ again is the phase function of $(4.27)$, by $(2.1)$ we have
$$|x^2p_y''(x)| \geq c |\lambda_1d_ix^{\alpha_i}| \eqno (4.43)$$
Then one can apply Lemma 4.2 for $n=2$ instead of $n=3$, and in place of $(4.39)$ we get a bound for $|T_i(\lambda_1,\mu_1,\mu_2)|$ of the form
$$ C' \|\phi\|_{C^1(V)}\bigg(\big|\{(x,y) \in D_i': x^{\alpha_i} < {1 \over |\lambda_1|}\}\big| + {1 \over |\lambda_1|^{1 \over 2}}
\int_{\{(x,y) \in D_i':\,\,x^{\alpha_i} \geq {1 \over |\lambda_1|}\}} {1 \over (x^{\alpha_i})^{1 \over 2}}\bigg)\eqno (4.44)$$
Then performing the argument analogous to before once again gives $(4.42)$, this time only using that $\epsilon < {1 \over 2}$.

Lastly, we consider  the case where $\epsilon = {1 \over 3}$ and  $\xi(x)$ does not have a zero of infinite order at $x = 0$.
If we had $|x^2p_y''(x)| \geq c|x^{\alpha_i}|$, then we could
proceed as in the case where $\xi(x)$ is identically zero since the argument required only that $\epsilon < {1 \over 2}$. However
this does not necessarily hold; this is because there can be an $x$ for which  $\lambda_1\alpha_i(\alpha_i - 1)d_ix^{\alpha_i} + \mu_2 s_i(s_i-1)\mu_2l_ix^{s_i} = 0$ and then the two main terms of $x^2p_y''(x)$ cancel.  However, since $s_i \neq \alpha_i$ there will always be an integer $j_0$ such that so long as $x$ is not in an interval $[2^{-j_0-1}, 2^{-j_0 + 1}]$ then one
does have $|x^2p_y''(x)| \geq c|x^{\alpha_i}|$. So we may apply the argument of the $\xi(x) = 0$ case for $T_{ij}$ with $j \neq j_0$ or $j_0 - 1$. Adding these in $j$ gives the bounds  $C \|\phi\|_{C^1(V)}|\lambda_1|^{-\epsilon}$  of $(4.42)$.

For $j = j_0$ or $j_0 - 1$, we apply the argument leading to $(4.37)$ unchanged. Then the steps leading to $(4.41)$ lead to the bound
$$ |T_{ij}(\lambda_1,\mu_1,\mu_2)| \leq C\|\phi\|_{C^1(V)}{1 \over |\lambda_1|^{1 \over 3}}\int_{\{t: t > {1 \over |\lambda_1|},\,\,2^{-j-1} < t < 2^{-j}\}} {1 \over 3} t^{-{4 \over 3}}t^{\epsilon} \,dt \eqno (4.45)$$
Since $\epsilon = {1 \over 3}$,  $(4.45)$ implies that we have
$$ T_{ij}|(\lambda_1,\mu_1,\mu_2)|  \leq C'\|\phi\|_{C^1(V)}{1 \over |\lambda_1|^{1 \over 3}} \eqno (4.46)$$
This is exactly the right-hand side of $(4.42)$. Since there are only two such $j$, if we combine with the earlier estimate for the the sum of the $T_{ij}$ for which $j \neq j_0$ or $j_0 - 1$, we see that $(4.42)$ once again holds and we are done.

\noindent This completes the proof of Theorem 1.2.

\noindent {\heading 5. Proofs of PDE Theorems.}

\noindent We start with the proof of Theorem 1.3.

\noindent {\bf Proof of Theorem 1.3.} 

\noindent Written on the Fourier transform side in the $x$ variables, $(1.6)$ becomes
$${\partial \hat{f} \over \partial t}(t,\xi_1,\xi_2) = iS(\xi_1,\xi_2) \hat{f}(t,\xi_1,\xi_2) $$
$$\hat{f}(0,\xi_1,\xi_2) = \hat{g}(\xi_1,\xi_2) \eqno (5.1)$$
This is solved by $\hat{f}(t,\xi_1,\xi_2) = e^{ itS(\xi_1,\xi_2)}\hat{g}(\xi_1,\xi_2)$. We are looking for solutions for when 
the support of $\hat{g}(\xi_1,\xi_2)$ is in a sufficiently small neighborhood of the origin. So we may fix a $\phi(\xi_1,\xi_2)$ 
supported in a neighborhood of the origin on which Theorems 1.1 and 1.2 hold such that $\phi(\xi_1,\xi_2) = 1$ on a neighborhood
$B$ of the origin, and we may assume that $\hat{g}(\xi_1,\xi_2)$ is supported on $B$. Thus we may write
$$\hat{f}(t,\xi_1,\xi_2)   = \hat{g}(\xi_1,\xi_2)e^{ itS(\xi_1,\xi_2)}\phi(\xi_1,\xi_2) \eqno (5.2)$$
Thus if $T(t,x_1,x_2)$ is as in Theorems 1.1 and 1.2 we have
$$f(t,x_1,x_2) = (g \ast T) (t,x_1,x_2)\eqno (5.3)$$
Here the convolution is in the $x$  variables for fixed $t$. By Theorems 1.1 and 1.2, one has 
$$|T(t,x_1,x_2)| \leq C\min\big((|t| + 2)^{-\epsilon}(\ln(|t| + 2))^m, |x|^{-{1 \over 2}}\big) \eqno (5.4)$$
 (We can add 2 to $|t|$ because $|T(t,x_1,x_2)|$ uniformly bounded
simply by taking absolute values of the integrand and then integrating.) In view of $(5.4)$, for a given $t$ it is natural to break up 
$T(t,x_1,x_2)$ into  $(|t| + 2)^{-\epsilon}(\ln(|t| + 2))^m < |x|^{-{1 \over 2}}$ and  $(|t| + 2)^{-\epsilon}(\ln(|t| + 2))^m \geq |x|^{-{1 \over 2}}$ pieces. To this end, for a given $t$ let $j_0$ be the nearest nonnegative integer to the $j$ for which $2^j = 
(|t| + 2)^{2\epsilon}(\ln(|t| + 2))^{-2m}$. As usual, let $\chi_{\{x: |x| < 2^{j_0}\}}(x)$ denote the characteristic function of the ball centered at the origin of radius $2^{j_0}$ and let $\chi_{\{x: 2^{j -1} \leq |x| < 2^{j}\}}(x)$ be the characteristic 
function of the annulus. Then we have
$$|T(t,x_1,x_2)| \leq C2^{-{j_0 \over 2}}\chi_{\{x: |x| < 2^{j_0}\}}(x_1,x_2) + C\sum_{j = j_0 + 1}^{\infty} 2^{-{j \over 2}}\chi_{\{x: 2^{j -1} \leq |x| < 2^{j}\}}(x_1,x_2) \eqno (5.5)$$
So by $(5.3)$ we have
$$|f(t,x_1,x_2)| \leq C 2^{-{j_0 \over 2}} \big||g|\ast \chi_{\{x: |x| < 2^{j_0}\}}(x_1,x_2)\big| + C\sum_{j = j_0 + 1}^{\infty}2^{-{j \over 2}}
\big||g| \ast \chi_{\{x: 2^{j -1} \leq |x| < 2^{j}\}}(x_1,x_2)\big| \eqno (5.6)$$
And therefore for any $q$, where the $L^q$ norm is in the $x$ variables we have
$$\|f(t,x_1,x_2)\|_q \leq C 2^{-{j_0 \over 2}} \||g| \ast \chi_{\{x: |x| < 2^{j_0}\}}(x)\|_q + C\sum_{j = j_0 + 1}^{\infty} 2^{-{j \over 2}}
\||g| \ast \chi_{\{x: 2^{j -1} \leq |x| < 2^{j}\}}(x)\|_q \eqno (5.7)$$
By Young's inequality, if ${1 \over q} = {1 \over p} + {1 \over r} - 1$, the above is bounded by
$$C 2^{-{j_0 \over 2} + {2j_0 \over r}} \|g\|_p + C\sum_{j = j_0 + 1}^{\infty} 2^{-{j \over 2} + {2j \over r} }\|g\|_p \eqno (5.8)$$
In the case that ${2 \over r} < {1 \over 2}$, the above converges and we get 
$$\|f(t,x_1,x_2)\|_q \leq C_q 2^{-{j_0 \over 2} + {2j_0 \over r}} \|g\|_p \eqno (5.9)$$
The condition that ${2 \over r} < {1 \over 2}$ translates into ${1 \over p} - {1 \over q} > {3 \over 4}$, which can occur when
$p < {4 \over 3}$. Given the definition of $j_0$, $(5.9)$ can be rewritten as 
$$\|f(t,x_1,x_2)\|_q \leq  C_q\big((|t| + 2)^{2\epsilon}(\ln(|t| + 2))^{-2m}\big)^{{2 \over r} - {1 \over 2}}\|g\|_p\eqno (5.10)$$
This is the same as
$$\|f(t,x_1,x_2)\|_q \leq  C_q(|t| + 2)^{4\epsilon({1 \over q} - {1 \over p} + {3 \over 4})}(\ln(|t| + 2))^{-4m({1 \over q} - {1 \over p} + {3 \over 4})}\|g\|_p\eqno (5.11)$$
This gives Theorem 1.3, except in the case where the exponent ${1 \over q} - {1 \over p} +  {3 \over 4}$ is zero. In this case, the exponent $-{j_0 \over 2} + 2{j_0 \over r}$ is always zero, so $(5.8)$ diverges. In this case, one may use the Hardy-Littlewood-Sobolev inequality instead,
directly using the bound $(5.5)$. For $(5.5)$ says that $|T(t,x_1,x_2)| \leq C|x|^{-{1 \over 2}}$, and since ${1 \over p} = {1 \over q} + {3 \over 4}$ the Hardy-Littlewood-Sobolev inequality gives that $\|f(t,x_1,x_2)\|_q  = \|g \ast T(t,x_1,x_2)\|_q \leq 
C_q\|g\|_p$, as long as $p \neq 1$ and $q \neq \infty$. This concludes the proof of Theorem 1.3.

\noindent In order to prove Theorems 1.4 and 1.5, in place of Theorems 1.1 and 1.2 we use the following, which generalizes 
estimates for the $\mu_1 = \mu_2 = 0$ case in [AGV].

\noindent {\bf Lemma 5.1.} Let $S(x,y)$ be as in Theorems 1.4 and 1.5, and define  $R(\lambda_1, \mu_1, \mu_2) $ by 
$$ R(\lambda_1, \mu_1, \mu_2)  = \int_{\R^2}e^{-\lambda_1S(x,y) +i\mu_1 x + i\mu_2y}\,\phi(x,y)\,dx\,dy \eqno (5.12)$$
There is a neighborhood of the origin $V$ that if the support of $\phi$ is contained in $V$ then for $\lambda_1 > 2$ one has the estimate
$$ |R(\lambda_1, \mu_1, \mu_2)| \leq C_S\|\phi\|_{C^1(V)} \min\big((\lambda_1)^{-\epsilon}(\ln \lambda_1)^m, |\mu|^{-1}\big) \eqno (5.13)$$ 
Here $(\epsilon,m)$ is as in $(4.2)$.

\noindent {\bf Proof.} We will first prove that $|R(\lambda_1, \mu_1, \mu_2)| \leq C_S\|\phi\|_{C^1(V)}(\lambda_1)^{-\epsilon}(\ln \lambda_1)^m$ and afterwards that $|R(\lambda_1, \mu_1, \mu_2)| \leq  C_S\|\phi\|_{C^1(V)} |\mu|^{-1}$. 
For the first estimate, we take absolute values on $(5.12)$ and integrate, obtaining
$$|R(\lambda_1, \mu_1, \mu_2)| \leq C_V\|\phi\|_{C^1(V)}\int_V e^{-\lambda_1S(x,y)}\,dx\,dy\eqno (5.14)$$
Note that 
$$\int_V e^{-\lambda_1S(x,y)}\,dx\,dy = \int_0^{\infty} \lambda_1 e^{-\lambda_1 t}|\{(x,y) \in V: S(x,y) < t\}\,dt \eqno (5.15)$$
We can truncate the integral here at $t = {1 \over 2}$ since the $t > {1 \over 2}$ portion gives an estimate much better than what we need.
Thus inserting $(4.2)$ we must bound 
$$\int_0^{1 \over 2} \lambda_1 e^{-\lambda_1 t} t^{\epsilon} |\ln t|^m\,dt \eqno (5.16)$$
Changing variables to $\lambda_1t$ here, the integral in $(5.15)$ is at most
$$\lambda_1^{-\epsilon}\int_0^{\lambda_1 \over 2} e^{-t}(|\ln t| + \ln \lambda_1)^m \eqno (5.17)$$
Regardless of whether $m = 0$ or $1$, equation $(5.17)$ is bounded by $C\lambda_1^{-\epsilon}(\ln \lambda_1)^m$. Putting
this back into $(5.14)$ gives the desired estimate.

We now prove that $|R(\lambda_1, \mu_1, \mu_2)| \leq  C_{S,\phi}  |\mu|^{-1}$. As before we can assume $|\mu| > 2$ as
the $|\mu| \leq 2$ case is obtained simply by taking absolute values and integrating. Rotating coordinates and shrinking our 
neighborhood $V$ of the origin if necessary, we assume that for some $k \geq 2$ we have $\partial_y^kS$ and $\partial_x^kS$ are
nonzero on $V$. Since the $x$
and $y$ axes are interchangeable here, without loss of generality we assume $|\mu_2| \geq |\mu_1|$.
We write the phase  in $(5.12)$ in the form
$$e^{-\lambda_1S(x,y) +i\mu_1 x + i\mu_2y} = (-\lambda_1\partial_yS(x,y) + i\mu_2)e^{-\lambda_1S(x,y) +i\mu_1 x + i\mu_2y} \times \bigg({1 \over -\lambda_1\partial_yS(x,y) + i\mu_2}\bigg) \eqno (5.18)$$
We integrate by parts in $(5.12)$, integrating the $(-\lambda_1\partial_yS(x,y) + i\mu_2)e^{-\lambda_1S(x,y) +i\mu_1 x + i\mu_2y}$ factor from $(5.18)$ and differentiating the rest. The derivative can land in two places. First, it can land on the $\phi(x,y)$ factor.
For this term, we take absolute values and integrate, using the bound $|{1 \over -\lambda_1\partial_yS(x,y) + i\mu_2}| \leq
{1 \over |\mu_2|} \leq {2 \over |\mu|}$, and we obtain the needed bound of $C\|\phi\|_{C^1(V)}|\mu|^{-1}$. The second
place the derivative can land is the ${1 \over -\lambda_1\partial_yS(x,y) + i\mu_2}$ factor, which becomes a factor of
${\lambda_1 \partial_{yy}S(x,y) \over  (-\lambda_1\partial_yS(x,y) + i\mu_2)^2}$. The resulting term is bounded in absolute value
by 
$$\int_{\R^2} {|\lambda_1 \partial_{yy}S(x,y)| \over |-\lambda_1\partial_yS(x,y) + i\mu_2|^2} |\phi(x,y)|\,dx\,dy \eqno (5.19)$$
We split this into two terms, depending on whether or not $|\lambda_1\partial_yS(x,y) | \geq |\mu_2|$. We get that $(5.19)$ is bounded by
$$\int_{|\lambda_1\partial_yS(x,y) | \geq |\mu_2|} {|\lambda_1 \partial_{yy}S(x,y)| \over (\lambda_1\partial_yS(x,y))^2 }|\phi(x,y)|\,dx\,dy$$
$$+\int_{|\lambda_1\partial_yS(x,y) | < |\mu_2|} {|\lambda_1 \partial_{yy}S(x,y)| \over \mu_2^2} |\phi(x,y)|\,dx\,dy
\eqno (5.20)$$
$$\leq C\|\phi\|_{C^1(V)}\int_{\{(x,y) \in V: |\lambda_1\partial_yS(x,y) | \geq |\mu_2|\}} {|\lambda_1 \partial_{yy}S(x,y)| \over (\lambda_1\partial_yS(x,y))^2 }\,dx\,dy$$
$$+ C\|\phi\|_{C^1(V)}{1 \over \mu_2^2} \int_{\{(x,y) \in V: |\lambda_1\partial_yS(x,y) | < |\mu_2|\}} |\lambda_1 \partial_{yy}S(x,y)| \,dx\,dy \eqno (5.21)$$
In the first term of $(5.21)$ we integrate in $y$ for fixed $x$. We use the fact $\partial_y^kS(x,y) \neq 0$ on $V$ to split the interval of
integration into boundedly many subintervals on which $\partial_{yy}S(x,y)$ has constant sign. Integrating the ${|\lambda_1 \partial_{yy}S(x,y)| \over (\lambda_1\partial_yS(x,y))^2 }$ leads to $\pm {1 \over \lambda_1\partial_yS(x,y)}$ at the boundary
points, which is bounded in absolute value by ${1 \over |\mu_2|}$ given that on the domain of integration one has 
$|\lambda_1\partial_yS(x,y) | \geq |\mu_2|$. Thus the overall term is bounded by $C\|\phi\|_{C^1(V)}{1 \over |\mu_2|} \leq 
C'\|\phi\|_{C^1(V)}{1 \over |\mu|}$ as needed.

In the second term of $(5.21)$, we do the analogous argument, and the resulting integration leads to $\pm \lambda_1\partial_yS(x,y)$ at the boundary
points of the subintervals of integration. This time the condition that $|\lambda_1\partial_yS(x,y) | < |\mu_2|$ on the domain of
integration leads to the term being bounded by  bounded by $C{1 \over \mu_2^2}\|\phi\|_{C^1(V)}|\mu_2| \leq  C'\|\phi\|_{C^1(V)}{1 \over |\mu|}$ and we are done.

\noindent {\bf Proof of Theorem 1.4.}

\noindent The proof will proceed much like the proof of Theorem 1.3. This time, on the Fourier transform
side the PDE becomes
$${\partial \hat{f} \over \partial t}(t,\xi_1,\xi_2) = -S(\xi_1,\xi_2) \hat{f}(t,\xi_1,\xi_2) $$
$$\hat{f}(0,\xi_1,\xi_2) = \hat{g}(\xi_1,\xi_2) \eqno (5.22)$$
So if $R(t,x_1,x_2)$ is as in Lemma 5.1, for $t \geq 0$ the equation is solved by
$$f(t,x_1,x_2) = (g \ast R) (t,x_1,x_2)\eqno (5.23)$$
By Lemma 5.1, we have
$$|R(t,x_1,x_2)| \leq C\min\big((t + 2)^{-\epsilon}(\ln(t + 2))^m, |x|^{-1}\big) \eqno (5.24)$$
This time we break into $T(t,x_1,x_2)$ into  $(t + 2)^{-\epsilon}(\ln(t + 2))^m < |x|^{-1}$ and  $(t + 2)^{-\epsilon}(\ln(t + 2))^m \geq |x|^{-1}$ pieces. We let $j_1$ be the nearest nonnegative integer to the $j$ for which $2^j = 
(t + 2)^{\epsilon}(\ln(t + 2))^{-m}$. Then in analogy to $(5.5)$ we have
$$|R(t,x_1,x_2)| \leq C2^{-j_1}\chi_{\{x: |x| < 2^{j_1}\}}(x_1,x_2) + C\sum_{j > j_1} 2^{-j_1}\chi_{\{x: 2^{j -1} \leq |x| < 2^{j}\}}(x_1,x_2) \eqno (5.25)$$
So in analogy to $(5.7)$, for any $q$ we have
$$\|f(t,x_1,x_2)\|_q \leq C 2^{-j_1} \||g| \ast \chi_{\{x: |x| < 2^{j_1}\}}(x)\|_q + C\sum_{j > j_1} 2^{-j_1}
\||g| \ast \chi_{\{x: 2^{j -1} \leq |x| < 2^{j}\}}(x)\|_q \eqno (5.26)$$
By Young's inequality, where ${1 \over q} = {1 \over p} + {1 \over r} - 1$, this is bounded by
$$C 2^{-j_1 + {2j_1 \over r}} \|g\|_p + C\sum_{j > j_1} 2^{-j + {2j \over r} }\|g\|_p \eqno (5.27)$$
So when $r > 2$ this converges and this time we have the bound
$$\|f(t,x_1,x_2)\|_q \leq C_q 2^{-j_1 + {2j_1 \over r}} \|g\|_p \eqno (5.28)$$
The condition that $r > 2$ translates into ${1 \over p} - {1 \over q} > {1 \over 2}$, which can occur when 
$p < 2$. By definition of $j_1$, $(5.28)$ is the same as 
$$\|f(t,x_1,x_2)\|_q \leq  C_q\big((t + 2)^{\epsilon}(\ln(t + 2))^{-m}\big)^{{2 \over r} - 1 }\|g\|_p\eqno (5.29)$$
This in turn is the same as
$$\|f(t,x_1,x_2)\|_q \leq  C_q\big((t + 2)^{2\epsilon}(\ln(t + 2))^{-2m}\big)^{{1 \over q} - {1 \over p} + {1 \over 2}}\|g\|_p\eqno (5.30)$$
This gives Theorem 1.4, except in the case where the exponent ${1 \over q} - {1 \over p} + {1 \over 2}$ is zero. In this case, $(5.25)$ gives that $|R(t,x_1,x_2)| \leq C|x|^{-1}$, and then the Hardy-Littlewood-Sobolev theorem gives that
$\|f(t,x_1,x_2)\|_q  = \|g \ast R(t,x_1,x_2)\|_q \leq C_{p,q,S}\|g\|_p$ as long as $p \neq 1$ and $q \neq \infty$  and we are done. 

\noindent {\bf Proof of Theorem 1.5.} 

\noindent On the Fourier transform side in the $x_1$ and $x_2$ variables, $(1.12)$ becomes
$$\hat{f}(\xi_1,\xi_2) = S(\xi_1,\xi_2)^{-\delta} \hat{g}(\xi_1,\xi_2) \eqno (5.31a)$$
Like in the previous two theorems, the support condition on $\hat{g}$ means we can insert a cutoff function in $(5.31a)$, turning the equation into
$$\hat{f}(\xi_1,\xi_2) =\hat{g}(\xi_1,\xi_2)  S(\xi_1,\xi_2)^{-\delta} \phi(\xi_1,\xi_2) \eqno (5.31b)$$
Thus if we define $Q(x_1,x_2)$ by
$$Q(x_1,x_2) = \int_{\R^2}  S(\xi_1,\xi_2)^{-\delta}e^{ix_1\xi_1 + ix_2\xi_2}\phi(\xi_1,\xi_2)\,d\xi_1\,d\xi_2 \eqno (5.32)$$
Then we have
$$f = g \ast Q \eqno (5.33)$$
Next, note that if $t > 0$ and $\delta > 0$ one has
$$\int_0^{\infty}e^{-ut}u^{\delta - 1}du = c_{\delta}t^{-\delta}\eqno (5.34)$$
Inserting this into $(5.32)$ gives
$$Q(x_1,x_2) = c_{\delta} \int_{\R^2} \bigg(\int_0^{\infty} e^{- uS(\xi_1,\xi_2)}u^{\delta - 1}\,du\bigg)e^{  ix_1\xi_1 + ix_2\xi_2}\phi(\xi_1,\xi_2)\,du\,d\xi_1\,d\xi_2$$
$$=  c_{\delta} \int_0^{\infty} u^{\delta - 1}\bigg(\int_{\R^2} e^{- uS(\xi_1,\xi_2) + ix_1\xi_1 + ix_2\xi_2}\phi(\xi_1,\xi_2)\,d\xi_1\,d\xi_2\bigg)du
\eqno (5.35)$$
We perform the $(\xi_1,\xi_2)$ integration in $(5.35)$, use the bounds from Lemma 5.1, then integrate the result in $u$. The result is
$$|Q(x_1,x_2)| \leq C \int_0^{\infty}u^{\delta-1}\min\big((u+ 2)^{-\epsilon}(\ln(u + 2))^m, |x|^{-1}\big)\,du \eqno (5.36)$$
We now bound $|Q(x_1,x_2)|$. For $|x| < 4$, we just use the bound obtained by taking absolute values and integrating in
$(5.32)$, and get $|Q(x_1,x_2)| < C$. Note that here we use that $\delta < \epsilon$; the fact that $(4.2)$ holds ensures that
$|S(\xi_1,\xi_2)|^{-{\delta}}$ is integrable.

Now assume $|x| > 4$. If $m = 0$, let $\epsilon' = \epsilon$, and if $m = 1$, let $\epsilon '$ be any number satisfying
$\delta < \epsilon' < \epsilon$. Then $(5.36)$ gives
$$|Q(x_1,x_2)| \leq C' \int_0^{\infty}u^{\delta-1}\min\big((u + 2)^{-\epsilon'}, |x|^{-1}\big)\,du \eqno (5.37)$$

$$= C'{1 \over |x|}\int_0^{|x|^{1 \over \epsilon'} - 2} u^{\delta - 1}\,du + C'\int_{|x|^{1 \over \epsilon'} - 2}^{\infty}
u^{\delta - 1}(u + 2)^{-\epsilon '} \,du \eqno (5.38)$$
The first term is bounded by  $C'{1 \over |x|}\int_0^{|x|^{1 \over \epsilon'}} u^{\delta - 1}\,du$, or $C''|x|^{{\delta \over \epsilon'} - 1}$. For the second term, we have
$$\int_{|x|^{1 \over \epsilon'} - 2}^{\infty}u^{\delta - 1}(u + 2)^{-\epsilon '} \,du \leq \int_{|x|^{1 \over \epsilon'} - 2}^{\infty}u^{\delta - \epsilon' - 1} du \eqno (5.39)$$
Given our assumption that $\delta < \epsilon'$, $(5.39)$ converges, and since $|x| > 4$, this is at most
$$C\int_{|x|^{1 \over \epsilon'}}^{\infty}u^{\delta - \epsilon' - 1} du \eqno (5.40)$$
This integrates to a term bounded by a constant times $|x|^{{\delta \over \epsilon '} - 1}$. Combining with the first term, we conclude that for $|x| \geq 4$ we have
$$|Q(x_1,x_2)| < C|x|^{{\delta \over \epsilon '} - 1} \eqno (5.41)$$
Combining with the $|x| < 4$, bound we have
$$|Q(x_1,x_2)| < C\min(1,|x|^{{\delta \over \epsilon '} - 1}) \eqno (5.42)$$
Note that the right-hand side of $(5.42)$ is in $L^r$ for $r > {2\epsilon' \over \epsilon' - \delta}$. Thus for such $r$, by Young's
inequality, if ${1 \over q} = {1 \over p } + {1 \over r} - 1$ one has a bound
$$\|f\|_q  = \|g \ast Q\|_q \leq C_{p,q,S}\|g\|_p \eqno (5.43)$$
One also has this bound when $r = {2\epsilon' \over \epsilon' - \delta}$ by the Hardy-Littlewood-Sobolev inequality, as
long as $p \neq 1$ and $q \neq \infty$. Stated in
terms of $p$ and $q$ alone, we have that an estimate of the form $(5.43)$ holds whenever ${1 \over q} \leq  {1 \over p} - {\delta \over 2 \epsilon'} - {1 \over 2}$, unless ${1 \over q} =  {1 \over p} - {\delta \over 2 \epsilon'} - {1\over 2}$ and $p = 1$ or
$q = \infty$. 

Given how $\epsilon'$ was defined, our conclusions are therefore as follows. When $m = 0$, there is an estimate of the form $(5.43)$ whenever  ${1 \over q} \leq  {1 \over p} - {\delta \over 2 \epsilon} - {1\over 2}$, except when ${1 \over q} =  {1 \over p} - {\delta \over 2 \epsilon} - {1\over 2}$ and $p = 1$ or $q = \infty$. When $m = 1$, there is an estimate of the form $(5.43)$ whenever  ${1 \over q} < {1 \over p} - {\delta \over 2 \epsilon} - {1\over 2}$. This concludes the proof of Theorem 1.5.

\noindent {\heading 6. References.}

\noindent [AGV] V. Arnold, S. Gusein-Zade, A. Varchenko, {\it Singularities of differentiable maps},
Volume II, Birkhauser, Basel, 1988. \parskip = 3pt\baselineskip = 3pt

\noindent [ArCuKa] G. I. Arhipov, V.N. Cubarikov, A.A. Karacuba, {\it Trigonometric integrals}, Izv. Akad. Nauk SSSR
Ser. Mat., {\bf 43} (1979), 971-1003, 1197 (Russian); English translation in Math. USSR-Izv., 15 (1980), 211-239.

\noindent [BNW] J. Bruna, A. Nagel, and S. Wainger, {\it Convex hypersurfaces and Fourier transforms},
Ann. of Math. (2) {\bf 127} no. 2, (1988), 333--365. 

\noindent [C] M. Christ, {\it Hilbert transforms along curves. I. Nilpotent groups}, Annals of 
Mathematics (2) {\bf 122} (1985), no.3, 575-596.

\noindent [CoMa] M. Cowling, Mauceri, {\it Oscillatory integrals and Fourier transforms of surface carried measures},
Trans. Amer. Math. Soc. {\bf 304} (1987), no. 1, 53-68. 

\noindent [D] J.J. Duistermaat, {\it Oscillatory integrals, Lagrange immersions, and unfolding of singularities}, Comm. Pure Appl.
Math., {\bf 27} (1974), 207-281.

\noindent [ESa] L. Erdos, M. Salmhofer, {\it Decay of the Fourier transform of surfaces with vanishing curvature}, (English summary) Math. Z. {\bf 257} (2007), no. 2, 261-294. 

\noindent [G1] M. Greenblatt, {\it The asymptotic behavior of degenerate oscillatory integrals in two
dimensions}, J. Funct. Anal. {\bf 257} (2009), no. 6, 1759-1798.

\noindent [G2] M. Greenblatt, {\it Resolution of singularities in two dimensions and the stability of integrals}, Adv. Math., 
{\bf 226} no. 2 (2011) 1772-1802.

\noindent [Gr] P. Gressman, {\it Uniform estimates for cubic oscillatory integrals},  Indiana Univ. Math. J. {\bf 57} (2008), no. 7, 3419-3442. 

\noindent [IKeM] I. Ikromov, M. Kempe, and D. M\"uller, {\it Estimates for maximal functions associated
to hypersurfaces in $\R^3$ and related problems of harmonic analysis}, Acta Math. {\bf 204} (2010), no. 2,
151-271.

\noindent [IM1] I. Ikromov, D. M\"uller, {\it On adapted coordinate systems}, Trans. AMS, {\bf 363} (2011), 2821-2848.

\noindent [IM2] I. Ikromov, D. M\"uller, {\it  Uniform estimates for the Fourier transform of surface-carried measures in
$\R^3$ and an application to Fourier restriction}, J. Fourier Anal. Appl, {\bf 17} (2011), no. 6, 1292-1332.

\noindent [K1] V. N. Karpushkin, {\it A theorem concerning uniform estimates of oscillatory integrals when
the phase is a function of two variables}, J. Soviet Math. {\bf 35} (1986), 2809-2826.

\noindent [K2] V. N. Karpushkin, {\it Uniform estimates of oscillatory integrals with parabolic or 
hyperbolic phases}, J. Soviet Math. {\bf 33} (1986), 1159-1188.

\noindent [PS] D. H. Phong, E. M. Stein, {\it The Newton polyhedron and
oscillatory integral operators}, Acta Mathematica {\bf 179} (1997), 107-152.

\noindent [S] E. Stein, {\it Harmonic analysis; real-variable methods, orthogonality, and oscillatory 
integrals}, Princeton Mathematics Series Vol. 43, Princeton University Press, Princeton, NJ, 1993.

\noindent [V] A. N. Varchenko, {\it Newton polyhedra and estimates of oscillatory integrals}, Functional 
Anal. Appl. {\bf 18} (1976), no. 3, 175-196.

\line{}
\line{}

\noindent Department of Mathematics, Statistics, and Computer Science \hfill \break
\noindent University of Illinois at Chicago \hfill \break
\noindent 322 Science and Engineering Offices \hfill \break
\noindent 851 S. Morgan Street \hfill \break
\noindent Chicago, IL 60607-7045 \hfill \break
\noindent greenbla@uic.edu \hfill\break

\end